\documentclass[10pt,twoside]{article}
\usepackage{latexsym,amsfonts,amssymb,stmaryrd,psfig}
\usepackage{graphicx}
\begin{document}

\newtheorem{de}{Definition}[section]
\newtheorem{thh}[de]{Theorem}
\newtheorem{lm}[de]{Lemma}
\newtheorem{pr}[de]{Proposition}
\newtheorem{co}[de]{Corollary}
\newtheorem{re}[de]{Remark}
\newtheorem{exa}[de]{Example}
\newtheorem{exe}[de]{Examples}
\newtheorem{ob}[de]{Remarks}

\renewcommand{\theequation}{\thesection.\arabic{equation}}

\newcommand{\enunt}[2]{\medskip\noindent\bf #1 \it #2 \rm \hfill $\Box$}
\newcommand{\te}[2]{\medskip\noindent\bf #1 \it #2 \rm}
\newcommand{\defi}[1]{\medskip\noindent\bf Definition \it #1 \rm}
\newcommand{\dem}[1]{\medskip\noindent\it Proof. \rm #1 \hfill $\Box$}
\newcommand{\demlung}[2]{\medskip\noindent\it #1 \rm #2 \hfill $\Box$}
\newcommand{\ex}[2]{\medskip\noindent\small\sc  #1\hspace{2mm}\rm #2 \normalsize\medskip}
\newcommand{\exh}[3]{\medskip\noindent\small\sc  #1\hspace{2mm}\rm #2 \\\noindent\footnotesize\sc Hint.\hspace{2mm}\rm #3 \normalsize\medskip}
\newcommand{\exhl}[7]{\medskip\noindent\small\sc  #1\hspace{2mm}\rm #2 \\\noindent\footnotesize\sc Hint.\hspace{2mm}\rm #3\\\medskip\noindent\sc #4 \hspace{2mm}\rm #5 \\\noindent\sc #6 \hspace{2mm}\rm #7 \normalsize\medskip}
\newcommand{\exs}[3]{\medskip\noindent\small\sc  #1\hspace{2mm}\rm #2 \\\noindent\footnotesize\sc Solution.\hspace{2mm}\rm #3 \normalsize\medskip}

\newcounter{nbib}
\newenvironment{bib}{\noindent\begin{list}{[\arabic{nbib}]}{\usecounter{nbib}\setlength{\parsep}{0mm}\setlength{\itemsep}{0mm}\setlength{\leftmargin}{7mm}\setlength{\rightmargin}{0mm}}}{\end{list}}

\newcounter{nlsn}
\newenvironment{lsn}[2]{\smallskip\noindent\begin{list}{\arabic{nlsn}}{\usecounter{nlsn}\setlength{\topsep}{0mm}\setlength{\itemsep}{#1}\setlength{\leftmargin}{#2}\setlength{\rightmargin}{\leftmargin}}}{\end{list}}

\newcounter{ntqftax}
\newenvironment{tqftax}[2]{\smallskip\noindent\begin{list}{\it(A\arabic{ntqftax})}{\usecounter{ntqftax}\setcounter{ntqftax}{0}\setlength{\topsep}{0mm}\setlength{\itemsep}{#1}\setlength{\leftmargin}{#2}\setlength{\rightmargin}{\leftmargin}}}{\end{list}}

\newcounter{nls}
\newenvironment{ls}[3]{\smallskip\noindent\begin{list}{#1}{\usecounter{nls}\setlength{\topsep}{0mm}\setlength{\itemsep}{#2}\setlength{\leftmargin}{#3}\setlength{\rightmargin}{\leftmargin}}}{\end{list}}


\def\too{\longrightarrow}
\def\Right{\Longrightarrow}
\def\rr{\rightrightarrows}
\def\ra{\rightarrarrow}
\def\hook{\hookrightarrow}
\def\dc{\Longleftrightarrow}
\def\lo{\longmapsto}

\def\T{\cal T\it}
\def\kk{\bf K\it}
\def\ll{\bf L\it}
\def\lll{\bf L\it}
\def\ttt{\bf T\it}
\def\oo{\bf O\it}
\def\cc{{\mathbb C}}
\def\zz{{\mathbb Z}}
\def\nn{{\mathbb N}}
\def\rr{{\mathbb R}}
\def\qq{{\mathbb Q}}
\def\dn{{\cal D}^n}
\def\an{{\cal A}^n}
\def\am{{\cal A}^m}
\def\aa{{\cal A}}
\def\a{{\mathfrak a}}
\def\dd{{\cal D}}
\def\fr{{\mathfrak R}}
\def\ahat{\hat{\cal A}}
\def\pfi{\varphi}
\def\gg{{\mathfrak g}}
\def\e{\{[e_i]\}_{i=1,\dots,2g+n-1}}
\def\f{\{[f_j]\}_{j=1,\dots,2g+n-1}}
\def\hs{H_1(S^3-F,\zz)}
\def\hf{H_1(F,\zz)}
\def\ao{\stackrel{\circ}{\aa(\emptyset)}}
\def\aoh{\widehat{\stackrel{\circ}{\aa(\emptyset)}}}
\def\agh{\widehat{{\cal A}(\G)}}
\def\cnh{\widehat{{\cal C}(n)}}
\def\on{{\cal O}_n}

\newcommand{\G}{\Gamma}
\newcommand{\cZ}{\check Z}
\newcommand{\Q}{{\mathbb Q}}
\newcommand{\bR}{{\mathbb R}}
\newcommand{\CP}{{\cal Q}}
\newcommand{\ZCP}{{\cal Q_{\mathbb Z}}}
\newcommand{\cP}{{\cal P}}
\newcommand{\cC}{{\cal E}}
\newcommand{\CC}{{\cal F}}
\newcommand{\ZCC}{{\cal F_{\mathbb Z}}}
\newcommand{\ZG}{{\cal G_{\mathbb Z}}}
\newcommand{\cA}{{\cal A}}
\newcommand{\hZ}{{\Hat Z}}
\newcommand{\ZA}{{\cA^{\mathbb Z}}}
\newcommand{\ZP}{{\cP^{\mathbb Z}}}
\newcommand{\cPe}{{\cal P}^e}
\newcommand{\cPo}{{\cal P}^o}
\newcommand{\cPec}{{\cal P}^{e,c}}
\newcommand{\cPoc}{{\cal P}^{o,c}}
\newcommand{\R}{{\mathbb R}}
\newcommand{\ve}{\varepsilon}
\newcommand{\Z}{{\mathbb Z}}
\newcommand{\cB}{{\cal B}}
\newcommand{\cS}{{\cal S}}
\newcommand{\sgn}{{\text sgn}}


\newcommand{\UP}{\raisebox{0pt}{
                 \begin{picture}(14,8)(-5,-4)
                 \thicklines
                 \put(0,2){\oval(16,16)[b]}
                 \end{picture}}}
\newcommand{\LP}{\raisebox{0pt}{
                 \begin{picture}(14,8)(-5,-4)
                 \thicklines
                 \put(0,-4){\oval(16,16)[t]}
                 \end{picture}}}
\newcommand{\U}{\raisebox{0pt}{
                 \begin{picture}(16,8)(-5,-4)
                 \thicklines
                 \put(0,2){\oval(16,16)[b]}    
                 \put(8,3){\vector(0,1){3}}
                 \end{picture}}}
\newcommand{\D}{\raisebox{0pt}{
                 \begin{picture}(16,8)(-5,-4)
                 \thicklines
                 \put(0,-3){\oval(16,16)[t]}    
                 \put(-8,-4){\vector(0,-1){3}}
                 \end{picture}}}
\newcommand{\arrowg}{\mathop{\raisebox{0pt}{
                 \begin{picture}(2,10)(0,-5)
                 \thicklines
                 \put(0,-5){\vector(0,1){10}}
                 \end{picture}}}_{[g]}}
\newcommand{\oGraph}{\raisebox{0pt}{
                 \begin{picture}(76,10)(-42,-5)
                 \thicklines
                 \put(-38,0){\circle{10}}
                 \put(-33,0){\vector(1,0){10}}
                 \put(-18,0){\circle{10}}
                 \put(-13,0){\vector(1,0){10}}
                 \put(-38,5){\vector(-1,0){3}}
                 \put(-18,5){\vector(-1,0){3}}
                 \put(-38,-5){\vector(1,0){3}}
                 \put(-18,-5){\vector(1,0){3}}
                 \dots
                 \put(0,0){\vector(1,0){10}}
                 \put(15,0){\circle{10}}
                 \put(15,5){\vector(-1,0){3}}
                 \put(15,-5){\vector(1,0){3}}
                  \end{picture}}}
\newcommand{\oOriented}{\raisebox{0pt}{
                 \begin{picture}(12,8)(-5,-4)
                 \thicklines
                 \put(0,0){\circle{10}}
                 \put(0,5){\vector(-1,0){3}}
                 \put(0,-5){\vector(1,0){3}}
                 \end{picture}}}

%
%
%
%

\title{\bf 3-cobordisms with their rational homology on the boundary
\thanks{The results of this article were obtained when the authors were at the Department of Mathematics, 
State University of New York at Buffalo, Buffalo, NY 14260-2900, USA}
}
\author{Dorin Cheptea and Thang T Q Le}
\date{}

\maketitle

\pagestyle{myheadings}
\markboth{Dorin Cheptea and Thang T Q Le}{3-cobordisms with their rational homology on the boundary}


\begin{abstract}
The object of this paper is to define a subcategory of the category of 3-cobordisms to which invariants of 
rational homology 3-spheres should generalize. We specify the notion of Topological Quantum Field Theory (in the sense 
of Atiyah [\ref{at88}]) to this case, and prove two intersesting properties that these TQFTs always have. In the case 
of the LMO invariant these properties amount to saying that the TQFT is anomaly-free.
\end{abstract}



\bigskip
\noindent 
This paper is organized as follows.
In Section 1 we define the topological categories of 3-cobordisms $\mathfrak Q, 
\mathfrak Z$, and $\mathfrak L$ that we want to consider, 
where the morphisms are certain {\it connected} oriented 3-cobordisms (with parametrizations 
of the boundary components). 
Their domain and range are also {\it connected} oriented closed surfaces.
All the proofs in this section are elementary when given rigurous definition.\footnote{Perhaps 
we should remove the proofs of Propositions 1 through 7 altogether.}
In Section 2 we give the axioms of TQFT for these categories, and prove two interesting properties.

There are two essential differences between this type of TQFTs and the classical
TQFTs, such as Reshetikhin-Turaev's for quantum invariants. Firstly,
we tailor the construction for integer and rational homology spheres. 
Therefore we restrict to {\it connected} cobordisms between {\it connected}
surfaces, and hence have no tensor product $\otimes$ structure induced
by disjoint union. Secondly, when gluing, we discriminate between 
the domain and the range of a cobordism. In particular, while we regard the
standard surface $\Sigma_g$ of genus $g$ in the domain as the boundary of the
standard handlebody $N_g$, we regard $\Sigma_g$ in the range as the boundary
of the complement of $N_g$ in $S^3$. Thus, gluing identically in our TQFT produces
$S^3$ as opposed to $\#_g(S^2\times S^1)$ in the Reshetikhin-Turaev TQFT.
This framework still allows to derive associated representations of the Torelli group,
in fact of a larger {\it Lagrangian subgroup} of the Mapping Class Group.

These definitions and properties are motivated by the construction of a
TQFT for the Le-Murakami-Ohtsuki invariant $Z^{LMO}$ [\ref{CL},\ref{MO}],
because this invariant is strong for rational homology three-spheres, and is 
weaker if the rank of homology is bigger. However the two simple properties that
we prove in section 2 can be useful for any TQFT on $\mathfrak Q$.

%
%
%
%

\section{The categories $\mathfrak Q\supset \mathfrak Z\supset \mathfrak L$ of semi-Lagrangian cobordisms}
\setcounter{equation}{0}

All maps and homeomorphisms in this paper are piecewise-linear, hence we will generally drop the term "PL".
Denote by $\Gamma^g$ and call {\em chain graph} (terminology borrowed from [\ref{MO}]) 
the abstract trivalent graph $\oGraph$ with oriented edges as indicated.
Label its subgraphs $\oOriented$ from $1$ to $g$ from left to right. 
For $g=1$, set $\Gamma^1=\oOriented$, one oriented edge, no vertices.
For $g=0$, set $\Gamma^0=$ one point.

\medskip
\defi{1) A pair $(\Gamma, R)$, consisting of an ordered disjoint union of chain graphs 
and an oriented surface with boundary, will be called a {\bf ribbon pair}
if it is the union of finitely many copies of the two pairs
depicted in the figure 1a, such that each point of $\Gamma$ has a
neighbourhood PL-homeomorphic to one of the two depicted pairs; the first set
in any pair is always a subset of the second; the union of the
first sets is  $\Gamma$; and the union of the second sets is $R$.
For $\Gamma=\Gamma^0$ we require $R$ to be homeomorphic to
$D^2$ and the point $\Gamma^0$ be the center of the disk.

2) A surface $R$ embedded in $M$, together with a subset
$\Gamma\subset R$, such that $(\Gamma, R)$ is a ribbon pair, and such that
at any trivalent vertex, viewed as a point in $M$, the three tangent vectors to $\Gamma$ 
coming from the three edges span a 2-plane,
shall be called a {\bf ribbon graph neighbourhood} of $\Gamma\subset M$, and $G=(\Gamma,R)$
will be called an {\bf embedded {\sc framed} graph}. 

3) Two {\bf triplets} $(K,G_1,G_2)$ and
$(L,H_1,H_2)$, consisting each of a framed oriented link $K$
(respectively $L$) in $S^3$, and two disjoint (and disjoint from the
corresponding link) embedded {\sc framed} graphs are {\bf equivalent} 
(notation $\cong$) if there is a PL-homeomorphism 
$\phi:S^3\rightarrow S^3$ which preserves the link and the embedded framed graphs, 
i.e. $\phi$ sends $K$ to $L$, the first embedded framed graph $G_1=(\Gamma_1, R_1)$ 
to the first embedded framed graph $H_1=(\Delta_1, S_1)$, and the second 
$G_2=(\Gamma_2, R_2)$ to the second one $H_2=(\Delta_2, S_2)$. Here 
$\emptyset$ is also considered a framed oriented link in $S^3$. Call 
$G_1=(\Gamma_1, R_1)$ the {\bf bottom}, and $G_2=(\Gamma_2, R_2)$ 
the {\bf top} of the triplet.}

\medskip
Fix $N_g$ -- a standard neighbourhood of $\Gamma^g$ in $S^3$.
$\Sigma_g=\partial N_g\subset S^3$ will be called the {\em standard oriented
surface of genus} $g$. Let $\overline{N_g}$ be the handlebody
complement of $N_g$ in $S^3$.  We also denote by
$\overline{\Gamma^g}$ the core of $\overline{N_g}$. Clearly
$\partial\overline{N_g}=-\Sigma_g$. When $g=0$, i.e. if
$\Gamma^g$ is a point, $N_g$ is a ball. Call $N_g$ {\em the
standard handlebody}, and $\overline{N_g}$ {\em the standard
anti-handlebody of genus $g$}.

Fix the $g$ pairs of standard loops $a_i$, $b_i$, $i=1,\dots,g$ on
$\Sigma_g$ as in figure 1b. Namely, $a_i$'s bound disks in  $N_g$,
the bounded component of $\rr^3 - \Sigma_g$,  while $b_i$'s bound
disks in $\overline{N_g}$. The embedded standard graph $\Gamma^g$
in $\rr^3\subset S^3$ has a preferred (the blackboard) ribbon
graph neighbourhood $R^g$ (with $\partial
R^g\subset\Sigma_g$). Hence we can construct its {\em push} along
the framing transversal to the ribbon, a graph
$\Gamma^g_{push}\subset\Sigma_g$ (together with a ribbon graph
neighbourhood $R^g_{push}\subset\Sigma_g$). Namely, the oriented
circle components of this graph coincide with the loops $b_i$.
Call $\Gamma^g_{push}$ {\em the standard b-graph}. Similarly, the
circle components of $\overline{\Gamma^g}_{push}$ are the loops $a_i$.
For "up-to-isotopy" pictures see figure 1. 

$\Sigma$ or $\Sigma_i$ will generically denote a standard $\Sigma_{g}$,
or a disjoint union of these.
For an oriented connected closed surface $S$, a {\em parametrization} 
of $S$ is an orientation-preserving homeomorphism between a 
standard $\Sigma$ and $S$.

\smallskip
\defi{Let $M$ be a compact oriented $3$-manifold with boundary
$\partial M = (-S_1)\cup S_2$, suppose also that parametrizations
$f_1$, $f_2$ of each $S_1$, $S_2$ are fixed. We will call such 
$(M,f_1,f_2)$ a (parametrized) {\bf (2+1)-cobordism}. 
$S_1$ will be called {\bf the bottom}, $S_2$ -- {\bf the top} of 
the cobordism. The cobordisms  $(M,f_1,f_2)$ 
and $(N,h_1,h_2)$ will be called {\bf equivalent (homeomorphic)} if there is a 
PL-homeomorphism $F:M\rightarrow N$ sending bottom to 
bottom and top to top preserving the parametrizations, 
i.e. $F\circ f_i=h_i$, $i=1,2$. We will use $\cong$ to 
denote equivalent cobordisms.}

\begin{figure}[tbp]
\centerline{
\psfig{file=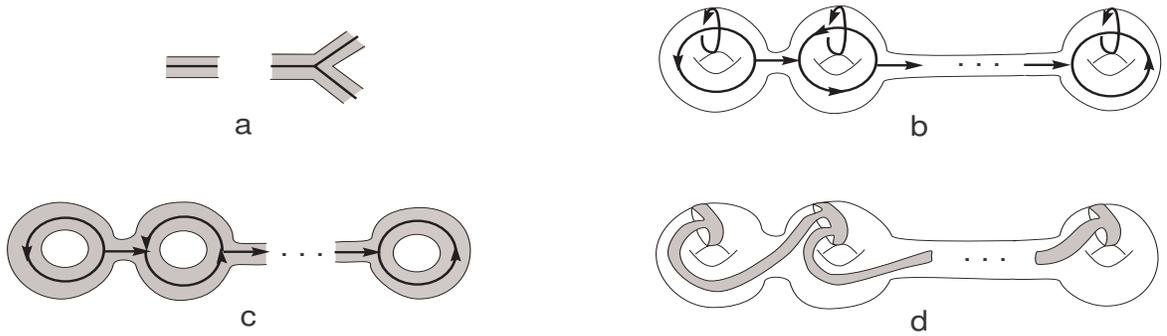,width=16cm,height=5cm,angle=0}
}
\caption{\sl {\bf a}: Local pieces that compose a ribbon pair,
{\bf b}: The standard b-graph, loops $a_i$, $b_i$ on $\Sigma_g$,
{\bf c}: $\Gamma^g\subset S^3$ and its preferred ribbon graph
neighbourhood $R^g$, {\bf d}:
$\overline{R^g_{push}}\subset\Sigma^g$, whose core (not drown
here) is the standard a-graph.} \label{fig1}
\end{figure}

\medskip
Using the parametrizations, we can glue the standard handlebody
$N_{g_1}$ to the bottom and the standard anti-handlebody
$\overline{N_{g_2}}$ to the top of $M$. Denote the result
$M\cup_{f_1}N_{g_1}\cup_{-f_2}(-\overline{N_{g_2}})$ by
$\widehat{M}$ and call it {\em the filling} of $(M,f_1,f_2)$. 


\subsection{Surgery description of gluing cobordisms}

Let $\mathfrak G$ denote set of equivalence classes of triplets 
$(L,G_1,G_2)$ in $S^3$. Let $\mathfrak C$ denote the 
set of equivalence classes of 3-cobordisms,
with non-empty bottom and top.

\medskip
\te{Proposition 1.}{1) There is a natural well-defined map $\kappa:{\mathfrak G}
\rightarrow\mathfrak C$ that associates to every equivalence class of triplets 
$(L,G_1,G_2)$ the equivalence class of cobordisms $(M,f_1,f_2)$, 
obtained by doing surgery on $L\subset S^3$, removing tubular
neighbourhoods $N_1,N_2$ of each $G_1,G_2$, and recording 
the parametrizations of the two obtained boundary components.  
If one glues according to these parametrizations a standard handlebody to 
$-\partial N_1$ and a standard anti-handlebody to $\partial N_2$,
then one obtains $S^3_L$.

2) $\kappa$ is surjective, and hence there exist maps $\upsilon:\mathfrak C\rightarrow 
\mathfrak G$ such that on 3-cobordisms with two 
non-empty boundary components $\kappa\circ\upsilon=id$.

3) Let a {\bf first Kirby move} on a triplet be the cancellation / insertion of a $\oo^{\pm 1}$ 
separated by an $S^2$ from anything else, and an {\bf extended (generalized) second Kirby move}
be a slide over a link component of an arc, either from another link component or from a chain graph.
Then if one factors $\mathfrak G$ by the extended Kirby moves and changes of orientations of link components, 
the induced map $\overline{\kappa}$ is a bijection.
} 

\medskip
\dem{1) The parametrizations are determined as follows. For $i=1,2$, 
let $N_i$ be the closure of a tubular neighbourhood of $G_i$ 
such that $\partial R_i\subset\partial N_i$. (Since $R_i$ is compact, 
there is a neighbourhood of $G_i$ of the form 
$R_i\times[-\varepsilon,\varepsilon]$. 
Take this as $N_i$.) Let $M=\overline{S^3_L-(N_1\cup N_2)}$, 
$S_1=\partial N_1$,  $S_2=-\partial N_2$. In particular, this 
construction yields $S_i\approx S^2$ if $G_i\approx\Gamma^0$. In the 
last case the possible parametrization is unique up to isotopy.

Fix a preferred point $x$ on each (open) upper half-circle of each
$\Gamma^g$. This determines a preferred point $x^{\prime}$ on each
circle component of $G_1$. The construction of $N_1$ produces
a preferred disk in $N_1$, centered at $x^{\prime}$, with boundary
in $\partial N_1$. Orient the boundary curves so that they twist
right-handedly with respect to the circle components of $G_1$. 
Similarly, fix a preferred point $y$ on each (open) lower half-circle of each 
$\Gamma^g$ and) construct an ordered system $b$
of oriented curves foor for the handlebody $N_2$.
(We assume $N_2$ contains the point at infinity $\infty\in S^3$,
hence we will refer to $N_2$ as the {\it anti-}handlebody.) 
Push $G_1$ along a framing transversal to $R_1$ to a graph on $S_1$, call it {\it the b-graph}. 
Analogously, push $G_2$ along a framing transversal $R_2$ to a graph on $S_2$, call it the 
{\it a-graph}.

If we cut-open $S_1$ along the b-graph and the system $a$, we get
an oriented surface homeomorphic to $D^2$. Define the
parametrization of $S_1$ by sending the standard b-graph on 
$\Sigma_{g_1}$ to the b-graph on $S_1$, and the system of loops 
$a_i$ on $\Sigma_{g_1}$ to the system $a$. Similarly for $S_2$: 
both $\Sigma-(\{standard\:\: a-graph\}\cup\{standard\:
\: system\:\: b\})$ and $S_2-(\{a-graph\}\cup\{system\:\: b\})$
are homeomorphic as {\it oriented} surfaces to the oriented
$D^2$. Any two orientation-preserving homeomorphisms of $D^2$ are
isotopic. Therefore the parametrization of $S_i$ is uniquely
determined up to isotopy by $G_i=(\Gamma_i, R_i)$. (Observe that
different choices of transversal framings to the ribbons lead to
isotopic b-graphs / a-graphs on $S_i$; different choices of the
points $x$ lead to isotopic systems $a$ / $b$.)

In conclusion, $\kappa$ of each triplet is well-defined as an
equivalence class of 3-cobordisms. It is obvious that 
via the above construction equivalent triplets yield the same
equivalence class of 3-cobordisms.

Moreover, in the above construction of parametrisations, the 
homeomorphism between the standard surface $\Sigma$ and $S$ extends 
to a homeomorphism between the standard handlebody / anti-handlebody 
and $N_i$. Hence $\widehat{M}\approx S^3_L$.

2) Let $\widehat{M}$ be the filling of a 3-cobordism
$(M,f_1,f_2)$, $\Gamma_1$, $\Gamma_2$ - the images in
$\widehat{M}$ of the cores and $R_1$, $R_2$ - the images in
$\widehat{M}$ of the preferred ribbon graph neighbourhoods of the
cores of the handlebody, respectively anti-handlebody of
$\widehat{M}$. Since $\widehat M$ is a closed 3-manifold, there is
a banded (unoriented) link $L\subset S^3$, such that $\widehat M\approx S^3_L$, the result of
surgery on $L$. Choose one such link $L$. Then there exist two
disjoint framed graphs $G_i$, $i=1,2$ in $S^3$, also disjoint 
from the link $L$, such that their remains after surgery coincide 
(up to ambient isotopy) with the pairs $(\Gamma_i, R_i)$, $i=1,2$ 
in $\widehat M$.\footnote{Indeed,
isotope the image of $(\Gamma_i, R_i)$ via $\widehat
M\stackrel{\approx}{\rightarrow}S^3_L$ to avoid the union of
surgery tori, a bounded nonseparating subset of $S^3_L$.}

3) see [\ref{MO}, proposition 2.1]
}

\begin{figure}[tbph]
\centerline{
\psfig{file=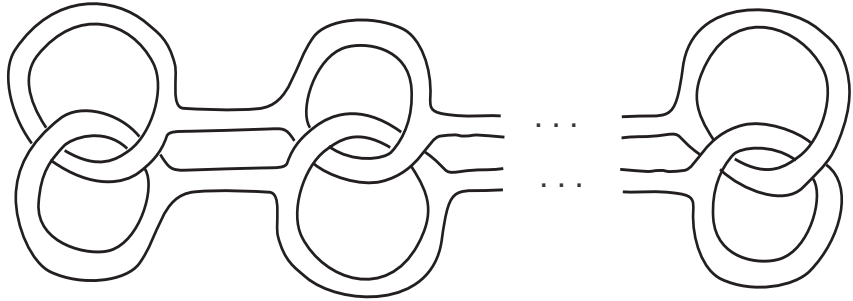,width=12cm,height=5cm,angle=0}
}
\caption{\sl The preferred choice of ribbons $R_i$, $i=1,2$ for $(\Sigma_g\times I,p_1,p_2)$.} \label{fig2}
\end{figure}

\medskip
To visualize the above proof it is helpful to imagine the standard b-graph and system $a$ 
(and respectively the standard a-graph and system $b$) cuting $\Sigma_g$ to a 2-disk. 
Strictly speaking the definitions of the b-graph and
system $a$ make sense if both $g_1,g_2>0$. However we can add the
remaining cases by making the following convention: if $\Gamma_i$
is a point, the b-graph/a-graph and the system $a$/$b$ are to be
considered the empty set. 
To draw a framed graph $G_i=(\Gamma_i, R_i)\subset S^3$, we only need to
draw the projection of $\Gamma_i$ on $\rr^2$,
which can be done in such a way that the preferred blackboard
framing determines $R_i$ up to isotopy.

Let us consider $\Sigma_g\times I\subset S^3 $. 
$\Sigma_g\times\{0\}$ and $\Sigma_g\times\{1\}$ are identified
with two very near (isotopic) copies of the standard embedding
$\Sigma_g\subset S^3$. $(\Sigma_g\times I, p_1, p_2)$ is a
3-cobordism, and its filling is homeomorphic to
$S^3$. The parametrizations of the bottom $p_1$ and top $p_2$ are
the ones induced via the isotopies in $S^3$ between the standard
$\Sigma_g$ and $\Sigma_g\times\{i\}$ from the identity
$id:\Sigma_g\rightarrow\Sigma_g\subset S^3$. 
To represent this 3-cobordism we choose framed graphs $R_1$, $R_2$ as
depicted in figure 2 
(projections on $\rr^2\subset\rr^3\subset S^3$). Observe the
consistency with Figure 1d.\label{p1p2}
Denote by $[D_{g}]_{L_0}$ the result of
surgery on the link $L_0$ in $D_{g}$, the three-dimensional manifold
depicted in figure 3. Note the disks $E_k$
and $E_k^{\prime}$ in $D_g$, that intersect $L_0$ once. If $g=0$,
set $D_{g}\approx S^2\times I$ and $L_0=\emptyset$.
$([D_{g}]_{L_0},p_1,p_2)$ is a 3-cobordism equivalent to
$(\Sigma_g\times I,p_1,p_2)$.

\begin{figure}[htb]
\centerline{
\psfig{file=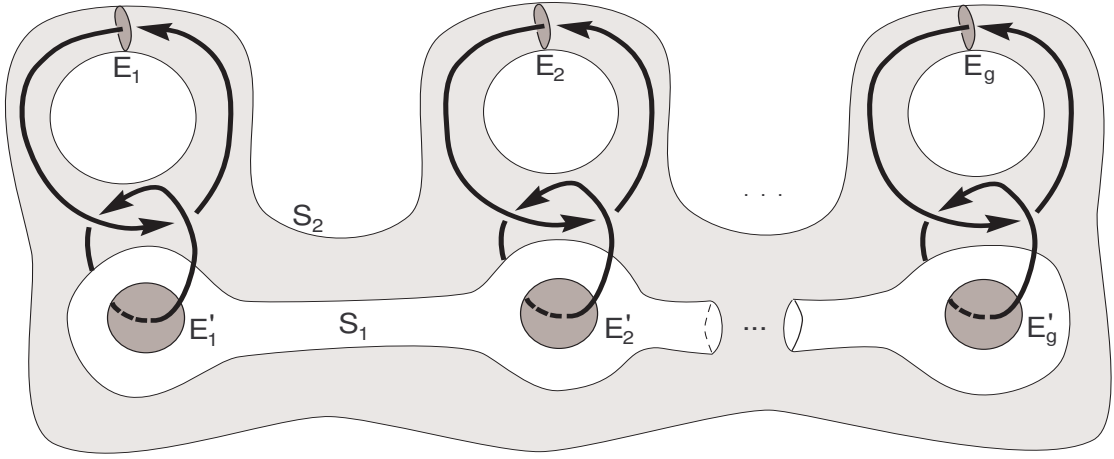,width=16cm,height=7cm,angle=0}
}
\caption{\sl $L_0\subset D_g$} \label{fig3}
\end{figure}

\medskip
\te{Proposition 2}{Let $(M, f_1, f_2)$ and $(N, h_1, h_2)$ be
arbitrary 3-cobordisms with connected bottoms and tops.
The following 3-cobordisms are equivalent:
\begin{equation}\label{CompCobordism}
(N\cup_{h_1\circ f_2^{-1}} M, f_1, h_2) \cong (N\cup_{h_1\circ
p_2^{-1}}[D_{g}]_{L_0}\cup_{p_1\circ f_2^{-1}} M, f_1, h_2)
\end{equation}
}

\vspace{-5mm}
\dem{Start with the right-hand-side. Using Kirby calculus, slide
the handles of the surface $S_1$ along the upper components of the
link $L_0$. Then the lower components of $L_0$ bound disks, so the
link can be canceled altogether. The two surfaces that remain are
clearly isotopic, and the parametrizations are equivalent since
both are images under isotopies of the identity parametrization of
$\Sigma_g\subset S^3$. The equivalence thus follows.}

\medskip 
Let $\Gamma$, respectively $\Gamma^{\prime}$ generically denote
the bottom, respectively the top of a triplet.
Call the union of the lower half-circles and the horizontal
segments of $\Gamma$, {\em the horizontal line} of $\Gamma$.
Similarly, call the union of the upper half-circles and the
horizontal segments of $\Gamma^{\prime}$, {\em the horizontal
line} of $\Gamma^{\prime}$. See figure 4a.

\begin{figure}[htb]
\centerline{
\psfig{file=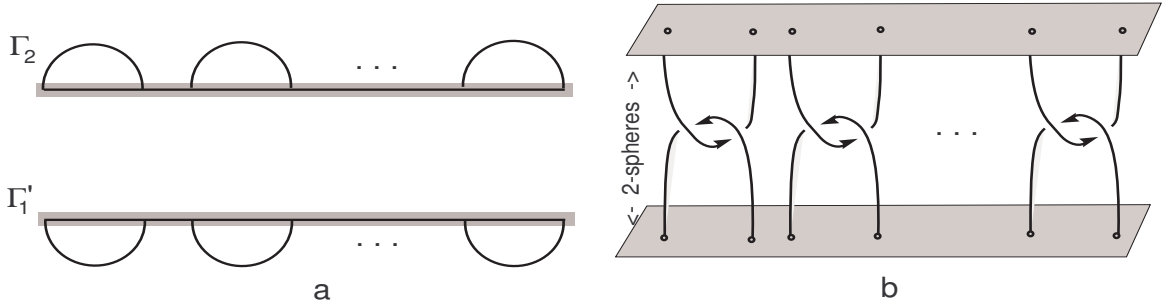,width=16cm,height=5cm,angle=0}
}
\caption{\sl {\bf a:} The horizontal segments of $\Gamma_1^{\prime}$ 
and $\Gamma_2$; {\bf b:} Framed tangle
$T_g\subset\overline{B(0,2)-B(0,1)}$.} \label{fig4}
\end{figure}

For a given genus $g$, let us decompose the manifold $D_g$,
together with the link $L_0$ inside, into a union of upper handles
$D^U$, lower handles $D^L$, and the rest of it $D^M$:
$D_{g}=D^U\cup D^M\cup D^L$. In figure 3 instead of the whole
handles, we have represented only disks $E_k$ and $E_k^{\prime}$,
whose neighbourhoods the handles are. Accordingly, $L_0\subset
D_{g}$ can be decomposed into three framed oriented tangles: an
oriented framed tangle in the upper handles, the oriented tangle
$T_g$ with the blackboard framing (see figure \ref{fig4}b), and an
oriented framed tangle in the lower handles. For $g=0$,
$T_0=\emptyset$.

\begin{figure}[htb]
\centerline{
\psfig{file=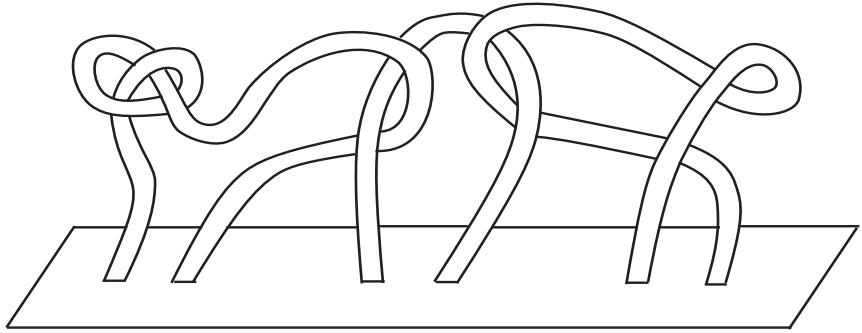,width=12cm,height=5cm,angle=0}
}
\caption{\sl An example (only the ribbons are depicted) of gluing
the upper handles of $D_g$.} \label{fig5}
\end{figure}

\medskip
\te{Proposition 3}{Let $(M_1,f_1,f_1^{\prime})$ and
$(M_2,f_2,f_2^{\prime})$ be two
3-cobordisms with connected non-empty bottoms and tops. Let
$\upsilon(M_1,f_1,f_1^{\prime})=(L_1,G_1,G_1^{\prime})$,
$\upsilon(M_2,f_2,f_2^{\prime})=(L_2,G_2,G_2^{\prime})$
where $\upsilon$ is some section of $\kappa$.
Remove a 3-ball neighbourhood of the horizontal line of
$G_1^{\prime}\subset S^3$, and identify the remain with 
$\overline{B(0,1)}$. Remove a 3-ball neighbourhood of the
horizontal line of $G_2\subset S^3$, and identify
the remain with $\overline{S^3-B(0,2)}$. Glue the framed tangle
$T_g\subset\overline{B(0,2)-B(0,1)}$ shown in figure \ref{fig4}b to
the ends of the remains of $G_1^{\prime}$ in $\overline{B(0,1)}$
and $G_2$ in $\overline{S^3-B(0,2)}$, strictly
preserving the order of the points, so that the composition of
these framed tangles is a smooth framed oriented link $L_0$ in
$S^3=(\overline{S^3-B(0,2)})\cup(\overline{B(0,2)-B(0,1)})
\cup(\overline{B(0,1)})$. Then 
\begin{equation}
\kappa(L_1\cup L_0\cup L_2, G_1, G_2^{\prime}) 
=  (M_2\cup_{f_2\circ(f_1^{\prime})^{-1}}M_1, f_1, f_2^{\prime})
\label{CompTriplet}
\end{equation}
where the ribbon neighbourhoods $R_1$, $R_2^{\prime}$ of $G_1$, $G_2^{\prime}$ 
in this formula are determined by the original $R_1$, 
$R_2^{\prime}$ in the two copies of $S^3$. Hence, any triplet
representing $(M_2\cup_{f_2\circ(f_1^{\prime})^{-1}}M_1, f_1, f_2^{\prime})$
is equivalent to $(L_1\cup L_0\cup L_2, G_1, G_2^{\prime})$ by
extended (generalized) Kirby moves and changes of orientations of link components.
}

\medskip
\dem{Let us first note that the result is obvious if the two
cobordisms are glued along a 2-sphere. For the general case, 
first use proposition 1.3 to "insert" $[D_g]_{L_0}$ ``between'' the two
cobordisms (for appropriate $g$). This changes
$(M_2\cup_{f_2\circ(f_1^{\prime})^{-1}}M_1, f_1, f_2^{\prime})$
to an equivalent 3-cobordism.

Observe that gluing the standard handlebody $N_g$ to $D_g$ along
$p_1$ produces a manifold homeomorphic to $N_g$. Moreover, gluing
$D_g\cup_{p_1}N_g$ along $f_2\circ p_2^{-1}$ to the
bottom of the cobordism $M_2$ produces a manifold
homeomorphic to the one obtained by gluing $N_g$ directly (along
$f_2$). In fact this homeomorphism is identity outside a
collar neighbourhood of the bottom of $M_2$.

Now, using the decomposition $D_g=D^U\cup D^M\cup D^L$, we note
that gluing $D_g\cup_{p_1}N_g$ to $M_2$ is the same as
first gluing $D^U$ on part of its boundary along the corresponding ``restriction'' of
the map $f_2\circ p_2^{-1}$, then gluing $(D^M\cup
D^L)\cup_{p_1}N_g\approx D^3$ along a 2-sphere.

Let us look at this glued 3-ball in the presentation of the filling
$\widehat{M_2}$ as $S^3_{L_2}$. By enlarging $D^U$
and $D^L$ if necessary, we may assume that inside this 3-ball there is
only a neighbourhood of the horizontal line of $\Gamma_2$
with $R_2$ of blackboard framing.

Apply a similar procedure to the top of the 3-cobordism
$M_1$: we may thus assume that $M_1\cup_{f_1^{\prime}\circ
p_1^{-1}}\overline{N_{g}}$ is equivalent (and the respective
homeomorphism is identity except in a collar neighbouhood of the
top of $M_1$) to a cobordism decomposed along a 2-sphere into
$M_1\cup_{restriction\: of\: f_1^{\prime}\circ p_1^{-1}}D^L$ and $D^M\cup
D^U\cup_{p_2^{-1}}\overline{N_{g}}$, the later homeomorphic to a
3-ball; and that the corresponding triplet in $S^3$ has
inside that 3-ball only a neighbouhood of the
horizontal line of $\Gamma_1^{\prime}$ with $R_1^{\prime}$ of blackboard framing.

Hence $(M_2\cup_{f_2\circ p_2^{-1}}[D_g]_{L_0}\cup_{p_1\circ(f_1^{\prime})^{-1}}M_1, f_1, f_2^{\prime})$
can be decomposed into three pieces: $M_2\cup_{restriction\:
of\: f_2\circ p_2^{-1}}D^U$, $D^M$, and
$D^L\cup_{restriction\: of\: p_1\circ(f_2^{\prime})^{-1}}M_1$, that have to be
glued along two 2-spheres. 
Observe that $D^M$ with the remaining tangle inside it is homeomorphic with
$\overline{B(0,2)-B(0,1)}$ with tangle $T_g$ inside.
Therefore (\ref{CompTriplet}), where all objects are as described in the
statement, holds.
}


\subsection{$\qq$- and $\zz$-cobordisms}

\noindent For a chain graph $\Gamma^g$ embedded in $S^3$,
denote $\mu_1,\dots,\mu_g$ the meridians of the upper
half-circles, defined by the right-hand rule, just as the
meridians of oriented link components are defined.

\medskip
\te{Proposition 4}{Let $\Gamma^g$ be embedded in an arbitrary
way in $S^3$, then $H_1(S^3-\Gamma^g,\zz)\cong\zz^g$, with free
generators $\mu_1,\dots,\mu_g$.}

\dem{Either write down the exact homology sequence of the pair
$(S^3, S^3-\Gamma^g)$  between $H_2(S^3,\zz)=0$ and
$H_1(S^3,\zz)=0$; or use Alexander duality; or write down the Mayer-Vietoris sequence (for
the reduced homology) of the decomposition $S^3-\Gamma^g=A\cup B$,
where $A$ is the complement in $D^3$ of the graph:

\medskip
\begin{center}
\begin{picture}(50,40)(90,-20)
\put(-50,0){\line(1,0){40}} \put(-50,30){\line(0,-1){30}}
\put(-40,0){\line(0,-1){30}} \put(-30,30){\line(0,-1){30}}
\put(-20,0){\line(0,-1){30}} \dots \put(10,0){\line(1,0){40}}
\put(20,30){\line(0,-1){30}} \put(30,0){\line(0,-1){30}}
\put(40,30){\line(0,-1){30}} \put(50,0){\line(0,-1){30}}
\put(150,0){$\left(\begin{array}{cccccc}\multicolumn{6}{c}{Id_{2g-1}}\\
1&-1&0&0&\ldots&0\\   0&0&1&-1&\ldots&0\\
\vdots&\vdots&\vdots&\vdots&\ddots&\vdots\\
0&0&0&0&\ldots&1\end{array}\right)$}
\end{picture}
\end{center}

\noindent
and $B$ is some embedding of $\sqcup_gI$ in $S^3-D^3$, to get
$H_1(S^3-\Gamma^g,\zz)$ as a quotient of $H_1(A,\zz)\oplus
H_1(B,\zz)\cong\zz^{2g-1}\oplus\zz^g$ through the image of
$H_1(A\cap B,\zz)\cong\zz^{2g-1}$ via the map given by the
matrix above.
}

\medskip
This proposition admits an obvious generalization to the case of
several connected components $\Gamma^{g_j}$. Loosely speaking, it
says that homology does not detect the horizontal lines.

\medskip
\te{Proposition 5} {Suppose $M$ is a connected compact oriented 
3-manifolds with two distinguished (not necessarily connected) 
boundary components $\partial M = (-S_1)\cup S_2$, let $f_1$, $f_2$ be
parametrizations of these surfaces, let $N_1$ and $\overline{N_2}$ be 
corresponding-to-the-genera disjoint unions of standard handlebodies, 
respectively anti-handlebodies, and let $i=(i_1,i_2): \partial
M\hookrightarrow M$ be the inclusion. The following conditions are
equivalent:

\medskip
\rm(1)\it \quad $ H_1(\widehat{M}, \zz)=0$

\rm(2)\it \quad $ H_1(M, \zz) = i_{\ast}(H_1(\partial M,\zz)/
(f_1,-f_2)_{\ast}H_1(N_1\sqcup -\overline{N_2},\zz))$ 

\medskip\noindent
They imply:

\medskip
\rm(3)\it \quad $ 2\cdot rank\: H_1(M;\zz) = rank\: H_1(\partial
M;\zz)$ }

\medskip
\dem{ We will prove this proposition for $S_1\approx\Sigma_{g_1}$ and
$S_2\approx\Sigma_{g_2}$, $g_1, g_2\geq 0$. The general case is
absolutely analogous.\footnote{The case when $S_j=\emptyset$
follows readily from the case $S_j\approx S^2$. Our convention then is that
$\emptyset_\ast=0$, where $\emptyset$ is the topological map and
$0$ is the algebraic one. Then $Im(i_j\circ f_j)_{\ast}=0$. In the case 
$S_1=S_2=\emptyset$, the condition $(2)$ reads $H_1(M,\zz)=0$.}

Applying Mayer-Vietoris to
$\widehat{M}=(M\cup_{f_1}N_{g_1})\cup_{-f_2}(-\overline{N_{g_2}})$,
using the fact that
$(M\cup_{f_1}N_{g_1})\cap\overline{N_{g_2}} = \partial
\overline{N_{g_2}}$ is connected, and then the second isomorphism
theorem, we obtain
$H_1(\widehat{M},\zz)\cong
\frac{H_1(M\cup_{f_1}N_{g_1},\zz)\oplus H_1(\overline{N_{g_2}},\zz)}
{i_{2\ast}H_1(S_2,\zz)}
\cong  \frac{H_1(M\cup_{f_1}N_{g_1},\zz)}{i_{2\ast}H_1(S_2,\zz)/
(i_2\circ(-f_2))_{\ast}H_1(\overline{N_{g_2}},\zz)}$. 
In a similar fashion $M\cap N_{g_1} = \partial N_{g_1}$ is connected,
and 
$H_1(M\cup_{f_1}N_{g_1},\zz)\cong\frac{H_1(M,\zz)}
{i_{1\ast}H_1(S_1,\zz)/(i_1\circ f_1)_{\ast}H_1(N_{g_1},\zz)}$.
Therefore
$H_1(\widehat{M},\zz)\cong\frac{H_1(M,\zz)}
{i_{1\ast}(H_1(S_1,\zz)/f_{1\ast}H_1(N_{g_1},\zz))
+i_{2\ast}(H_1(S_2,\zz)/(-f_2)_{\ast}H_1(\overline{N_{g_2}},\zz))}$, 
which proves $(1)\Longleftrightarrow (2)$.

$(1)\Right(3)$ We will give a geometric proof, naturally extending 
linking relations from the case of {\it closed} 3-manifolds [\ref{GS}].

There is a link $L\subset S^3$ such that $\widehat
M\approx S^3_L$. Moreover, this link can be taken disjoint from
the embedding $\Gamma^{g_1}\sqcup\overline{\Gamma^{g_2}}\hookrightarrow S^3$,
such that (identifying as before $\Gamma^{g_1}\sqcup\overline{\Gamma^{g_2}}$
with the remain after surgery on $L$) $\widehat
M-N(\Gamma^{g_1}\sqcup\overline{\Gamma^{g_2}})\approx M$ for a tubular
neighbourhood $N(\Gamma^{g_1}\sqcup\overline{\Gamma^{g_2}})$. Hence, there is
a neighbourhood $N(L\sqcup\Gamma^{g_1}\sqcup\overline{\Gamma^{g_2}})$ such
that $M$ is obtained from
$S^3-N(L\sqcup\Gamma^{g_1}\sqcup\overline{\Gamma^{g_2}})$ by adding a
2-handle and a 3-handle for each component $L$. By proposition 4
(and the remark afterwords)
$H_1(S^3-N(L\sqcup\Gamma^{g_1}\sqcup\overline{\Gamma^{g_2}}),
\zz)\cong\zz^{|L|+g_1+g_2}$ with generators
$\mu_1,\dots,\mu_{|L|+g_1+g_2}$, the right-handed oriented
meridians of the link components $K_j$, the upper-half-circles
$U_j$ of $\Gamma^{g_1}$, and the lower half-circles $V_j$ of
$\overline{\Gamma^{g_2}}$. Therefore $H_1(M,\zz)$ is a quotient of the former
by $|L|$ relations, one for each 2-handle. If the component $K_i$
of $L$ has surgery coefficient (framing) $l_{i}$ and the preferred
longitude $\lambda_i$ (i.e $\lambda_i$ has framing $0$ in $S^3$),
then the corresponding relation is $l_{i}\mu_i+\lambda_i=0$. But
$\lambda_i$ is the boundary of a Seifert surface $F_i$ in $S^3$,
punctured by the other components $K_j$, $j\not=i$, $U_j$ and
$V_j$ (clearly $F_i$ can be taken disjoint from two 3-balls, hence
the intersection with $\Gamma^{g_1}-U_j$ and $\Gamma^{g_1}-V_j$
can be avoided). These intersections result in a surface
$F_i^{\prime}$ with additional boundary components homologous to
$\pm\mu_j$.
Therefore $F_i$ determines a relation $\lambda_i
=\sum\limits_{j\not=i}lk(K_j,K_i)\mu_j
+\sum\limits_{j}lk(U_j,K_i)\mu_{|L|+j}
+\sum\limits_{j}lk(V_j,K_i)\mu_{|L|+g_1+j}$ for the homology of
$M$. Hence $H_1(M,\zz)$ is a quotient of $\zz^{|L|+g_1+g_2}$
through the image of $A:\zz^{|L|}\rightarrow\zz^{|L|+g_1+g_2}$
given by the $(|L|+g_1+g_2)\times|L|$-matrix 
$\left(\begin{array}{c}lk_{ij}\\ {\bf lk_{\Gamma L}}\end{array}\right)$. 
(Here $lk_{ii}=l_{i}$.)

On the other side, adding a 2-handle (along the corresponding
$\mu_j$) for each component $U_j$ and $V_j$, as well as a 3-handle
for $\Gamma^{g_1}$ and 3-handle for $\Gamma^{g_2}$, one obtains
$\widehat M$. At the level of homology this adds precisely the
relations $\mu_j=0$ for $j=|L|+1,\dots,|L|+g_1+g_2$. This means
that $H_1(\widehat M,\zz)$ is a quotient of $\zz^{|L|}$ through
the image of the linear homomorphism $B$ given in the $\mu$-basis
by the $|L|\times|L|$-matrix $(lk_{ij})$, the linking matrix of
$L$. Since $H_1(\widehat M,\zz)=0$, $(lk_{ij})$ is unimodular,
hence invertible (over $\zz$). Therefore there is a basis
$(\nu_1,\dots,\nu_{|L|})$ of the module freely generated by
$\mu_1,\dots,\mu_{|L|}$, such that the matrix representing $B$ in
the new basis has the form $\left(\begin{array}{ccc} \pm 1&&{\bf 0}\\
&\ddots&\\ {\bf 0}&&\pm 1\end{array}\right)$. Correspondingly
$A=B\oplus\bf lk_{\Gamma L}$ has in the new basis the form
$\left(\begin{array}{ccc}\pm 1&&{\bf 0}\\ &\ddots&\\ {\bf 0}&&\pm 1\\
&{\bf lk^{\prime}_{\Gamma L}}&\end{array}\right)$. 
Therefore, when writing down the relations in $H_1(M,\zz)$ for the 
system of generators 
$\nu_1,\dots,\nu_{|L|},\mu_{|L|+1},\dots,\mu_{|L|+g_1+g_2}$, the 
generators  $\nu_j$, $j=1,\dots,|L|$ can be eliminated, together with all relations,
without adding any new relations. Hence
$H_1(M,\zz)\cong\zz^{g_1+g_2}$, freely generated by
$\mu_{|L|+1},\dots,\mu_{|L|+g_1+g_2}$. The statement $(3)$ now
follows.
}

\medskip 
\te{Proposition 6} {Suppose $M$ is a connected compact oriented 
3-manifolds with two distinguished (not necessarily connected) 
boundary components $\partial M = (-S_1)\cup S_2$, let $f_1$, $f_2$ be
parametrizations of these surfaces, let $N_1$ and $\overline{N_2}$ be 
corresponding-to-the-genera disjoint unions of standard handlebodies, 
respectively anti-handlebodies, and let $i=(i_1,i_2): \partial
M\hookrightarrow M$ be the inclusion. The following conditions are
equivalent:

\medskip
\rm(1)\it \quad $ H_1(\widehat{M}, \qq)=0$

\rm(2)\it \quad $ H_1(M, \qq) = i_{\ast}(H_1(\partial M,\qq)/
(f_1,-f_2)_{\ast}H_1(N_1\sqcup -\overline{N_2},\qq))$

\medskip\noindent
They imply:

\medskip
\rm(3)\it \quad $ 2\cdot rank\: H_1(M;\qq) = rank\: H_1(\partial
M;\qq)$ }

\medskip
\dem{is identical to the proof of the previous proposition in
every aspect, except that in the course of proving $(1)\Right(3)$,
$H_1(\widehat M,\qq)=0$ implies only that the matrix of $B$ is
invertible (over $\qq$). The basis $\nu_1,\dots,\nu_{|L|}$ is then
over $\qq$, i.e. it is a linear combination of
$\mu_1,\dots,\mu_{|L|}$, but in general only with rational
coefficients. Again, these $\nu$-generators can be eliminated
together with all relations.}

\medskip
Of cause, the conditions (3) in propositions 5 and 6 are equivalent.

A 3-cobordism satisfying the equivalent conditions (1), (2) of Proposition 5 (respectively 6)
will be called a {\it $\zz$-cobordism} (respectively a {\it $\qq$-cobordism}). 
In both definitions we have allowed one
or both $S_i$ to be empty, although from the point of TQFT the case of empty
top and/or bottom is indistinguished from the case when that component is $S^2$.

\medskip
\te{Corollary 7}{If $H_1(\widehat{M},\zz)=0$, then $H_1(M,\zz)$
is free of rank sum of the genera of its boundary components. More 
generally, $H_1(M,\zz)$ can have only the kind of torsion 
$H_1(\widehat{M},\zz)$ has.}

\medskip
\dem{The proof of proposition 6 can be repeated for $\zz_p$. Hence
$rank\:H_1(M,\zz_p)=rank\:H_1(M,\qq)$ for all $p$ for which 
$H_1(\widehat{M},\zz_p)=0$, which implies that then $p$-torsion can not 
occur. In particular if $M$ is a $\zz$-cobordism, this is true for all $p$.
Apply the structure theorem for finitely generated abelian groups.
}

\medskip
In general $(3)$ does not imply $(1)$ in the statement of the
above propositions. For example, let $M$ be the manifold obtained 
from $S^3$ by excising one component of the Hopf link in $S^3$ and 
performing surgery on the other. Then $(3)$ is true, while $(1)$
is not. But, if we restrict to 3-cobordisms those equivalence class 
is in the category $\mathfrak L$ below, then such
phenomena are excluded a priori (Proposition 9).

For connected $\partial M$ condition $(2)$ clearly implies $H_{\ast}(M; \partial M)=0$.
However, for example $\kappa(\emptyset,G,G^{\prime})$ where $G$ are $G^{\prime}$ are
the components of the two-component unlink in $S^3$, obviously satisfies $(1)$,
and hence $(2)$, but fails to satisfy $H_{\ast}(M; S_i)=0$, the second condition from
the definition of an h-cobordism [\ref{milnor}].


\subsection{Description of the categories}

We will be interested in three categories, $\mathfrak Q\supset
\mathfrak Z\supset\mathfrak L$,
which we now describe.  
Objects in each of these are {\em natural numbers}. The
morphisms between $g_1$ and $g_2$ are certain equivalence
(homeomorphism) classes of connected 3-cobordisms
with connected bottom $S_1$ of genus $g_1$ and connected top $S_2$ of genus $g_2$. The
composition-morphism is the equivalence class of the 3-cobordisms in
(\ref{CompCobordism}).
The equivalence classes of the cobordisms
$([D_g]_{L_0},p_1,p_2)$, $g\geq 0$ play the r\^ole of identity in these 
categories. 

Let us first describe some additional objects. Let
$L^a$, $L^b$ be the submodules of $H_1(\Sigma_g, \zz)$ generated by
$a_i$'s and $b_i$'s respectively. Each is a Lagrangian submodule
with respect to the algebraic intersection form $\omega$ on
$H_1(\Sigma_g, \zz)$. Suppose $M$ is a $\qq$-cobordism,
with boundary $(-S_1)\cup S_2$. Let $L_i$, $i=1,2$ be submodules of
$H_1(S_i,\zz)$, viewed inside $H_1(M,\zz)$. We say that 
$L_1 \ge L_2$ in $M$, if every element in $L_2$ is 
$\zz$-homologically equivalent to an element in $L_1$. (Or, intuitively, $L_2$, 
which is in the homology of the component $S_2$, can be ``moved'' 
to the part $L_1$ of the homology of the component $S_1$.) 
Similarly we define $L_1\leq L_2$. If $L_1\geq L_2$ and $L_1\leq L_2$, 
we have $L_1=L_2$. Using the parametrizations of the boundary 
components, we can speak about the submodules $L^a_i$, $L^b_i$ of 
$H_1(S_i)$, $i=1,2$ , which correspond to $L^a$, $L^b$.

The morphisms in $\mathfrak Z$ from $g_1$ to $g_2$
are equivalence classes of $\zz$-cobordisms with boundary $(-S_1)\cup
S_2$ such that $g(S_i)=g_i$ and
\begin{equation}\label{LagrangeCondition}
L^a_1 \geq L^a_2,\; {\rm and} \; L^b_1 \leq L^b_2
\end{equation}

\noindent Such $\zz$-{\bf cobordisms} will be called {\bf semi-Lagrangian}.
An absolutely similar construction 
works for $\mathfrak Q$, just replace $\zz$ by $\qq$. 

\begin{figure}[t]
\centerline{
\psfig{file=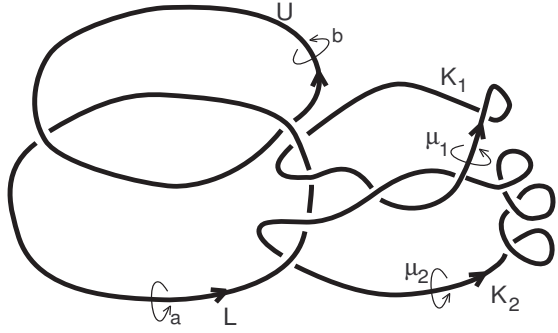,width=9cm,height=5cm,angle=0}
}
\caption{\sl An example to illustrate that condition (\ref{LagrangeCondition})
over $\qq$ is weaker than over $\zz$.} \label{fig6}
\end{figure}

\medskip
\te{Example}{In general condition (\ref{LagrangeCondition})
over $\zz$ is stronger than (\ref{LagrangeCondition}) over $\qq$.}
Let $(M,f_1,f_2)$ be a representative of the equivalence class 
of 3-cobordisms obtained by applying $\kappa$ to the 
triplet $(K=K_1\cup K_2,L,U)$ shown in figure \ref{fig6}. Let 
$a,b,\mu_1,\mu_2$ be the corresponding meridians. Then $H_1(M,\zz)$
is the abelian group with the classes of these as generators, and 
relations $\mu_1+\mu_2+a=0$ and $-3\mu_2+\mu_1-a=0$. They imply that
$2(\mu_1-\mu_2)=0$, i.e. $\mu_1-\mu_2$ is a torsion element. 
$L_1^a,L_1^b,L_2^a,L_2^b$ are generated respectively by 1 element each, 
the classes of respectively $a,L,U,b$. Note that $L$ is homologous to 
$\mu_1-\mu_2-b$, hence over $\zz$, $L_1^b\not\leq L_2^b$, but over $\qq$,
$\mu_1-\mu_2=0$, hence $L_1^b=L_2^b$. In this example 
$H_1(\widehat{M},\zz)\cong\zz_4$, i.e. $(M,f_1,f_2)$ is a $\qq$-cobordism, 
but not a $\zz$-cobordism. {\it If one assumes, however, that $(M,f_1,f_2)$ 
is a $\zz$-cobordism, then conditions (\ref{LagrangeCondition}) over $\zz$ and 
$\qq$ are equivalent.} This follows from (2) in proposition 5, and 
from the fact that for $\zz$-cobordisms, $L_2^a$ and $L_1^b$ are free (by 
Proposition 7).

\medskip
Note that if condition (\ref{LagrangeCondition}) is satisfied with 
equalities, then necessarily $g_1=g_2$. ( (\ref{LagrangeCondition}) is
defined with the assumption that $M$ is already a $\qq$-cobordism.)

\medskip
\te{Example}{Condition (\ref{LagrangeCondition}) may hold with strict inclusion.}
Consider $\kappa(\emptyset,G,G^{\prime})$, where $G$ is the braid-closure of a
generator of the braid group $B_2$, and $G^{\prime}$ is the braid axis (in $S^3$).
In the case of integer homology both condition in (\ref{LagrangeCondition}) are
strict, as it is easy to check. This example also shows that {\it 3-cobordisms
representing elements of the category $\mathfrak Z$ do not necessarily satisfy
$H_{\ast}(M; S_i)=0$.} (They are both $\cong\zz_2$ in this example.) It is not hard
to see that if a homology-cobordism triad $(M,S_1,S_2)$ is enhanced with
parametrizations such that we get a $\qq$-cobordism, then the second condition from the definition 
of an h-cobordism [\ref{milnor}] implies (\ref{LagrangeCondition}), and in fact with equalities. Hence,
{\it if we restrict $\qq$-cobordism, requiring (\ref{LagrangeCondition}) contains all
homology-cobordisms, and more.}

\medskip
\te{Proposition 8}{The composition of two morphisms (say, class of $M$ and class of $N$) in 
category $\mathfrak Q$ (resp. $\mathfrak Z$) is again a 
morphism in the category $\mathfrak Q$ (resp. $\mathfrak Z$), i.e. $\mathfrak Q$ and $\mathfrak Z$ are categories.}

\medskip
\dem{The fact that both $M, N$ represent $\qq$-cobordism (respectively $\zz$-cobordism)
means (by propositions 5 and 6) that all 1-dimensional
homology (over $\mathbb Q$, respectively over $\zz$) in $M$ can be
considered as in the boundary. So when we glue $M$ with $N$ along
a surface $S$, the 1-dimensional homology is either in the top
component of $N$, in the bottom component of $M$, or in the
``middle'' surface $S$. But now, by condition
(\ref{LagrangeCondition}), all the cycles of type $L^a$ in $S$ can
be moved down to the bottom, and all the cycles of type $L^b$ can be
moved up to the top. So the 1-dimensional homology is still
sitting on the boundary. Since for any 3-cobordism 
$(M,f_1,f_2)$, $(i_1\circ f_1)_{\ast}H_1(N_{g_1},\zz)$ is $L_1^b$, 
and $(i_2\circ f_2)_{\ast}H_1(\overline{N_{g_2}},\zz)$ is $L_2^a$,
using (\ref{LagrangeCondition}) one can verify condition (2) in
Proposition 6 (respectively 5). Therefore 
$N\cup_{f_N\circ(f_M^{\prime})^{-1}}M$ is a $\qq$-cobordism (respectively 
a $\zz$-cobordism). For $N\cup_{f_N\circ(f_M^{\prime})^{-1}}M$ one can 
easily check (\ref{LagrangeCondition}).
}

\medskip
The morphisms in the category $\mathfrak L$ we define to be the
equivalence classes of 3-cobordisms of the form
$(M=\Sigma_g\times I,f,f^{\prime})$, where $f,f^{\prime}\in
Aut(\Sigma_g)$, which are in $\mathfrak Q$. The following 
proposition shows that 
$\mathfrak Q\cap \mathfrak L = \mathfrak Z\cap \mathfrak L$, 
hence to require the equivalence class of $(\Sigma_g\times[0,1],f,f^{\prime})$
to be a morphism  in category $\mathfrak Q$ or $\mathfrak Z$ 
is equivalent. Clearly, in $\mathfrak L$ there are no morphisms 
between non-equal natural numbers. 
Recall that we use the following notation for the indices 1,2 and 
$^{\prime}$ : $(M_1,f_1,f_1^{\prime})(M_2,f_2,f_2^{\prime})=
(M_2\cup_{f_2\circ(f_1^{\prime})^{-1}}M_1,f_1,f_2^{\prime})$.

\medskip
\te{Proposition 9}{Consider a 3-cobordism
$M=(\Sigma_g\times[0,1],f\times 0,f^{\prime}\times 1)$, where the
parametrization of the top differs by that of the bottom by the
automorphism $w=(f^{\prime})^{-1}\circ f$. Then its equivalence
class depends only on the isotopy class of $w$ (i.e. we don't need
to specify both $f$, $f^{\prime}$), and the following are
equivalent:

\medskip
\rm(1)\it \quad the equivalence class of $(\Sigma_g\times[0,1],f\times
0,f^{\prime}\times 1)$ is  a morphism  in the category $\mathfrak
L$

\rm(2)\it \quad $L^a = w_*(L^a)$ and $L^b = w_*(L^b)$,

\medskip
\noindent In particular, $\widehat{M}$ is a $\zz$-homology sphere.
}

\medskip
\dem{$(\Sigma_{g}\times I, f\times 0, f^{\prime}\times 1)$ is
equivalent as a 3-cobordism to $(\Sigma_{g}\times I,
((f^{\prime})^{-1}\circ f)\times 0=w\times0, id \times 1)$,
therefore $\widehat{M}\cong\widehat{(\Sigma_g\times I,w\times
0,id\times 1)}$. As a 3-manifold, $\widehat{M}$ is homeomorphic
to $\overline{N_g}\cup_w N_g$. Using
Mayer-Vietoris theorem for the decomposition $\overline{N_g}\cup_w
N_g$ of $\widehat{M}$,  we can see that condition (2)
already ensures that $M$ is a $\zz$-cobordism. Hence $(2)\Right(1)$.

$(1)\Right(2)$. Since $w$ is a homeomorphism, $w_*$ is an automorphism,
hence by (\ref{LagrangeCondition}) $w_*$ has to be an automorphism over 
$\qq$ of each $L^a$, $L^b$. 
On the other side $w_{\ast}(L^a)\subset H_1(M,\zz)=L^a\oplus L^b$, 
and as a $\zz$-submodule of a free module, $w_{\ast}(L^a)$ must be free. 
Using this and the fact that $w_{\ast}(L^a)=L^a$ over $\qq$, we conclude
$w_{\ast}(L^a)=L^a$ over $\zz$. 
Similarly $w_{\ast}(L^b)=L^b$ over $\zz$.
}

\medskip
Suppose a closed 3-manifold is the result of gluing a 
standard handlebody $N_g$ to the standard anti-handlebody 
$\overline{N_g}$ along a homeomorphism $w$ of the standard surface 
$\Sigma_g$, whose action on the homology (in the $a_1,\dots,a_g,b_1,\dots,b_g$
basis) is given by a symplectic matrix 
$\left(\begin{array}{cc}A&B\\ C&D \end{array}\right)$.
By Mayer-Vietoris theorem this 3-manifold is a $\zz$HS if and only if
$A$ is invertible. But the set of symplectic matrices 
$\left(\begin{array}{cc}A&B\\ C&D \end{array}\right)$ with $A$ invertible
is not closed under multiplication, as it is easy to notice from the 
following example: $\left(\begin{array}{cc}1&1\\ 0&1 \end{array}\right)
\left(\begin{array}{cc}1&0\\ -1&1 \end{array}\right)=
\left(\begin{array}{cc}0&1\\ -1&1 \end{array}\right)$. The result of gluing
along each of the first two homeomorphisms (of $S^1\times S^1$)  is $S^3$, 
while the result of gluing along their composition is $S^1\times S^2$.

\medskip
\enunt{Corollary 10}{ The composition of two cobordisms
$(\Sigma_{g}\times I, f_1\times 0, f_1^{\prime}\times
1)\cong(\Sigma_{g}\times I, w_1\times 0, id\times1)$ and
$(\Sigma_{g}\times I, f_2\times 0, f_2^{\prime}\times
1)\cong(\Sigma_{g}\times I, w_2\times 0, id\times1)$ along
$(f_2\times 0)\circ((f_1^{\prime})^{-1}\times 1)$ (respectively
$(w_2\times 0)\circ(id\times 1)^{-1}$) is the 3-cobordism
$(\Sigma_{g}\times I,f_2\circ(f_1^{\prime})^{-1}\circ f_1\times
0,f_2^{\prime}\times 1) \cong (\Sigma_{g}\times
I,(f_2^{\prime})^{-1}\circ f_2\circ(f_1^{\prime})^{-1}\circ
f_1\times 0,id\times 1) \cong (\Sigma_{g}\times I, (w_2 \circ
w_1)\times 0, id\times 1)$. In particular, the composition of two
morphisms of category $\mathfrak L$ is again a morphism in the
same category.
}

\medskip
\defi{Denote by ${\cal L}_g$ the subgroup of the Mapping Class Group,
consisting of isotopy classes of elements $w\in Aut(\Sigma_g)$
such that $w_*(L^a) = L^a$ and $w_*(L^b) = L^b$ (over $\qq$ or
over $\zz$, is equivalent by the previous proposition), and call it
the {\bf Lagrangian subgroup of the MCG}.
}

\medskip
TQFTs based on $\mathfrak Q$ induce representations of 
${\cal L}_g$. This subgroup of $MCG(g)$ is big enough to be 
interesting, it contains the Torelli group. In fact its image under
the action on homology is the group of matrices of the form
$\left(\begin{array}{cc}A&0\\ 0&(A^T)^{-1}\end{array}\right)$,
where $A\in GL(g,\zz)$. This subgroup of $Sp(2g,\zz)$ is not normal.

Proposition 9 above also shows that within this particular type
($M=\Sigma_g\times I$ and $w\in{\cal L}_g$), statement (3) from
propositions 5 and 6 implies statement (1).

\medskip\noindent
{\it Remark.} Let $\lambda$ denote the Casson invariant of homology 3-spheres. By fixing the standard handlebody 
of genus $g$ in $\rr^3\subset S^3$ we fixed a Heegaard homeomorphism that Morita [\ref{morita}] calls $\iota_g$, 
and by taking the filling $\widehat{(\Sigma_g\times I,\varphi,id)}$ we obtain a manifold denoted by Morita $W_{\pfi}$.
Every $\pfi$ is a composition of Dehn twists $\tau_{\gamma}^{\pm 1}$.
Using Proposition 2, we can "insert" $[D_g]_{L_0}$ between every two twists, or put it another way, express
$(\Sigma_g\times I,\varphi,id)$ as a composition cobordisms $(\Sigma_g\times I,\tau_{\gamma}^{\pm 1},id)$. 
Every twist can be replaced with $\pm 1$-surgery on a knot $K_i$. Hence we obtain $W_{\pfi}$ as surgery
on the link $L=K_1\cup\dots\cup K_n\cup L_{0_1}\cup\dots\cup L_{0_{n-1}}$, such that if one removes $L_0$'s, then
the remaining $K_1\cup\dots\cup K_n$ is split.
If $\pfi\in{\cal K}_g$, the kernel of the Johnson homomorphisms,
then it is a composition of Dehn twists $\tau_{\gamma}^{\pm 1}$ on separating simple closed
curves,\footnote{The problem of finding an algorithm to express an arbitrary $\pfi\in{\cal K}_g$ as a
product of $\tau_{\gamma}^{\pm 1}$'s is unsolved.}
i.e. $\gamma$ has zero linking number with every circle component of $L_0$'s.

\begin{figure}[thbp]
\centerline{
\psfig{file=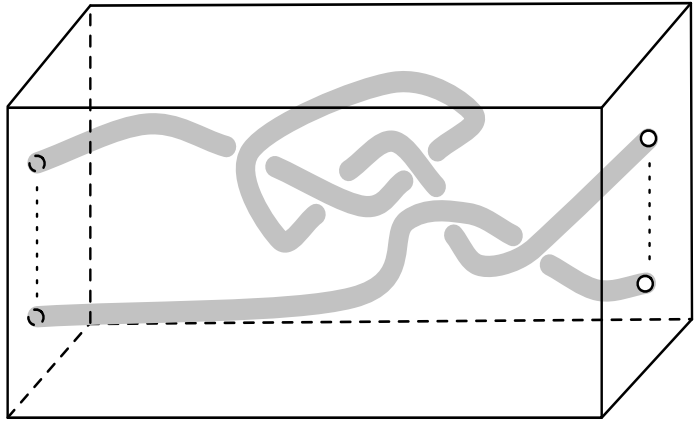,width=10cm,height=7cm,angle=0}
}
\caption{\sl $(L,G)$ less a tubular neighbourhood of the horizontal line of $G$.} \label{multiplicationBox}
\end{figure}

\noindent
{\it Remark.}
Let ${\mathfrak C}_{\emptyset}$ denote the set
of connected 3-cobordisms with empty bottom
and connected top.
For $(M_1, \emptyset, f_1),
(M_2, \emptyset, f_2)\in {\mathfrak C}_{\emptyset}$
and $(L_1,G_1)\subset S^3$, $(L_2,G_2)\subset S^3$ 
such that $\kappa(L_1,G_1)=(M_1,\emptyset,f_1)$,
$\kappa(L_2,G_2)=(M_2,\emptyset,f_2)$\footnote{A pair $(L,G)$ is 
similar to a triplet $(\emptyset,L,G)$.},
we can remove a tubular neighbourhood of the horizontal line of
$G_1$ ($=$ a ball in $S^3$), respectively of $G_2$,
and glue the two "boxes" as in figure 7, from left to right, 
afterwards filling back in the standard way a horizontal line.
Denote the result
by $(L_1\cup L_2,G_1\bullet G_2)$, and define:

\begin{equation}\label{multiplication}
(M_1, \emptyset, f_1)\bullet(M_2, \emptyset, f_2)=\kappa(L_1\cup L_2,G_1\bullet G_2)
\end{equation}

\noindent
Observe that the new 3-cobordism does not depend on the 
choice of pairs $(L_i,G_i)$.
In the case of $g=0$, $\bullet$ is the connected sum, i.e. this 
operation is another way (alternative to composition of cobordisms) of generalizing the connected sum.
Note that $\widehat{(M_1,\emptyset,f_1)\bullet(M_2,\emptyset,f_2)}=
\widehat{(M_1,\emptyset,f_1)}\#\widehat{(M_2,\emptyset,f_2)}$.
In particular the sets ${\mathfrak C}_{\emptyset}\cap\{\zz{\rm -cobordism}\}$
and  ${\mathfrak C}_{\emptyset}\cap\{\qq{\rm -cobordism}\}$ are closed under $\bullet$.
However, the Kontsevich-LMO invariant of 3-cobordisms [\ref{CL},\ref{MO}] is not multiplicative 
with respect to $\bullet$, except in some paricular cases.

\medskip\noindent
{\it Remark.}
One can consider a slightly modified version of the category $\mathfrak Q$, where the objects are connected
parametrized surfaces (instead of natural numbers) and morphisms are certain classical (i.e. not parametrized)
3-cobordisms. Our choice here was motivated by our main example of TQFT (LMO), where $\mathfrak Q$ gives the
simpliest formulations. In the next section one can replace it by a "$\mathfrak Q$-like category" instead.

%
%

\section{TQFTs for the categories $\mathfrak Q\supset\mathfrak Z\supset\mathfrak L$}
\setcounter{equation}{0}

Define a TQFT $(\T,\tau)$ based on the cobordism category 
$\mathfrak Q$ (or a subcategory of it, or a $\mathfrak Q$-like category) to be 1) a covariant functor $\T$ 
from the category those objects are the objects of ${\mathfrak Q}$ 
(i.e. natural numbers) and morphisms are the 
homeomorphisms of parametrized surfaces to a subcategory
${\mathfrak V}_K$ of the category of $K$-modules, such that 
$\T(0)=K$, where $K$ is a commutative module with a conjugation operation; 
and 2) a map $\tau$ that associates
to each 3-cobordism $(M,f_1,f_2)$ a $K$-homomorphism 
$\tau(M):\T(\Sigma_1)\rightarrow\T(\Sigma_2)$,
satisfying the following axioms:

\it
\medskip
\begin{tqftax}{0mm}{10mm}

\item (Naturality) If $(M_1,\Sigma_1,\Sigma_1^{\prime})$, 
$(M_2,\Sigma_2,\Sigma_2^{\prime})$ are two 3-cobordisms, 
and $f:M_1\rightarrow M_2$ is a homeomorphism of 3-cobordisms, 
preserving the parametrizations, then the following diagram is 
commutative:
$$\T(\Sigma_1)\stackrel{\tau(M_1)}{\longrightarrow}
\T(\Sigma_1^{\prime})$$
$$\T(f|_{\Sigma_1})\downarrow\hspace{2cm}
\downarrow\T(f|_{\Sigma_2})$$
$$\T(\Sigma_2)\stackrel{\tau(M_2)}{\longrightarrow}
\T(\Sigma_2^{\prime})$$

\item (Functoriality) If $M_1$, $M_2$ are 3-cobordisms, 
$f=f_2\circ(f_1^{\prime})^{-1}:\partial_{top}(M_1)\rightarrow
\partial_{bottom}(M_2)$ is the gluing homeomorphism, and denote 
$M=M_2\cup_f M_1$, then $\tau(M)=k\cdot\tau(M_2)\circ\tau(M_1)$.
$k\in K$ is called {\sc the anomaly}.

\item (Normalization) Let 
$(\Sigma\times[0,1],(\Sigma\times0,p_1),(\Sigma\times1,p_2))$ 
be the  3-cobordism mentioned on page \pageref{p1p2}, then  
$$\tau(\Sigma_g\times[0,1],(\Sigma_g\times0,p_1),
(\Sigma_g\times1,p_2))=id_{\T(\Sigma_g)}$$

\item (pseudo-Hermitian structure) There is a superstructure on each 
element $V$ of ${\mathfrak V}_K$, i.e. it admits an antimorphism 
$\overline{\cdot}:V\rightarrow V$ (a map linear in 0-supergrading 
and antilinear in 1-supergrading), that commutes with ($=$ is natural 
with respect to) surface homeomorphisms. There is a canonical map $V\rightarrow V^{\ast}$,
which composed with the above antimorphism extends (from the particular 
case when $\T(\Sigma_1)=K$) to an antimorphism 
$\overline{\cdot}:Mor(\T(\Sigma_1),\T(\Sigma_2))\rightarrow
Mor(\T(-\Sigma_2),\T(-\Sigma_1))$, that commutes with 
homeomorphisms of 3-cobordisms, such that
$$\tau(-M)=\overline{\tau(M)}$$
\end{tqftax}

\rm
\noindent
We can not require multiplicativity or self-duality since in 
the category $\mathfrak Q$ all cobordisms are connected. 
Conditions (A1-A3) say that $\tau:{\mathfrak Q}\rightarrow\cal A$ is a
pseudo-functor. $\tau$ would is a true functor when there is no anomaly.
If the set of $\tau(M,S^2,\Sigma)$'s, spans (in the closure for infinite-dimensional
modules) $\T(\Sigma)$, the  TQFT is called {\it non-degenerate}.


\subsection{The behavior of the signature and the determinant}

\medskip
\te{Proposition 11}{Let $(M_1,f_1,f_1^{\prime})$ and 
$(M_2,f_2,f_2^{\prime})$ be two  3-cobordisms. 
Suppose $(M_1, f_1, f_1^{\prime}) = \kappa(L_1,G_1,G_1^{\prime})$,
$(M_2,f_2,f_2^{\prime}) = \kappa(L_2,G_2,G_2^{\prime})$, 
and $(M_2\cup_{f_2\circ(f_1^{\prime})^{-1}} M_1, f_1,
f_2^{\prime}) = \kappa(L_1\cup L_0\cup L_2, G_1,G_2^{\prime})$,
the later triplet obtained from the previous two by the construction
described in Proposition 3.
Denote $\sigma^1_+=sign_+(lk(L_1))$, $\sigma^2_+=sign_+(lk(L_2))$,
$\sigma_+=sign_+(lk(L_1\cup L_0\cup L_2))$, and let $g$ be the
genus of the connected closed surface along which is this splitting.
Then the integer $s(M,M_1,M_2):=\sigma^1_++\sigma^2_++g-\sigma_+$ is an invariant of the decomposition
$M=M_2\cup_{f_2\circ(f_1^{\prime})^{-1}} M_1$, i.e. it does not depend on the choice of triplets 
representing the 3-cobordisms $M_1$ and $M_2$.
}

\medskip
\dem{Suppose we have another such choice of triplets.
The new unoriented links $L_1^{\prime}$ and $L_2^{\prime}$ are related to 
$L_1$ and $L_2$ by a finite sequence of Kirby moves and changes of orientations
of link components. Each such move can be 
thought of also as a move from $L_1^{\prime}\cup L_0\cup L_2^{\prime}$
to $L_1\cup L_0\cup L_2$. If a K-1 move changes $\sigma_+^1$ by $\pm 1$, then so
does it to $\sigma_+$, and hence $\sigma_+^1+\sigma_+^2+g-\sigma_+$ remains 
unchanged. Similarly happens if a K-1 move changes $\sigma_+^2$. 
K-2 moves do not affect the signature. Indeed, recall that the operations 
of the type ``add a $j^{th}$ row and column to an $i^{th}$ row an column'' 
(which corresponds to sliding $i^{th}$ component over the $j^{th}$ component of 
the link), do not change the signature of a symmetric matrix. They correspond to 
re-writing a symmetric bilinear form in a different basis, and Sylvester's inertia 
theorem applies. If the orientation of a link component of $L_1$ or $L_2$ changes,
it corresponds to multiplying both linking matrices $lk(L_i)$ and $lk(L_1\cup L_0\cup L_2)$
on the left and on the right by a diagonal matrix with all entries $1$ except one entry $-1$.
Hence again by Sylvester's theorem all signatures involved are invariant.
}

\medskip
Note that $s(M,M_1,M_2)$ is not an invariant of 
a triad (in the sense of Milnor [\ref{milnor}], i.e. if one ignores parametrizations. 

Let $(L,G,G^{\prime})$ be a triplet and $(M,f,f^{\prime})=\kappa(L,G,G^{\prime})$. 
Recall that (see Proposition 5) we can talk about linking number between a link component $K$ 
and a circle $U$ of a chain graph, as well as between two circles $U$ and $V$ of chain graphs:
$lk(K,U)=lk(U,K)=$ the linking number between $K$ and the knot obtained 
from the graph by deleting all but the circle component $U$; and similarly for $lk(U,V)$.
We can then define {\it the linking matrix of a triplet}:
\begin{equation}\label{linkingMatrixTriplet}
lk(L,G,G^{\prime})=\left(\begin{array}{ccc}
lk(L) & lk(L,G) & lk(L,G^{\prime}) \\
lk(G,L) & lk(G,G) & lk(G,G^{\prime}) \\
lk(G^{\prime},L) & lk(G^{\prime},G) & lk(G^{\prime},G^{\prime})
\end{array}\right)=\left(\begin{array}{ccc}
A & B^T & C^T \\
B & D & E^T \\
C & E & F
\end{array}\right)
\end{equation}

\noindent
where $A, D, F$ are symmetric matrices.

Let $\mu$ be the column-vector consisting of the meridians of $L$, $m$ the column-vector of
the meridians of $G$, and $m^{\prime}$ be the column-vector of the meridians of $G^{\prime}$. 
Then $H_1(M,\zz)$ is the $\zz$-module generated by the elements of $\mu, m, m^{\prime}$ 
with $|L|$ relations -- the elements of the column-vector $A\mu+B^Tm+C^Tm^{\prime}$.
The semi-Lagrangian conditions $L_1^b\leq L_2^b$, $L_2^a\leq L_1^a$ can be expressed
$\zz\!<\!B\mu+Dm+E^Tm^{\prime}\!> \:\leq \zz\!<\!m^{\prime}\!>$, respectively 
$\zz\!<\!C\mu+Em+Fm^{\prime}\!> \:\leq \zz\!<\!m\!>$.
If $\widehat{M}$ is a $\qq$- or $\zz$-homology sphere, $A$ is invertible over $\qq$. Moreover,
$A^{-1}\in{\cal M}_{|L|\times|L|}(\zz\left[\frac{1}{\det A}\right])$. 
Therefore $\mu=-a^{-1}B^Tm-A^{-1}C^Tm^{\prime}$. (This equality of column-vectors with entries in
$\zz\left[\frac{1}{\det A}\right]$ has to be read over $\zz$, i.e. to mean that multiplying an entry on 
the left with the denominator of the corresponding entry on the right gives the numerator on the 
right side.) Hence the semi-Lagrangian conditions can be expressed
$\zz\!<\!(D-BA^{-1}B^T)m+(E^T-BA^{-1}C^T)m^{\prime}\!> \:\leq \zz\!<\!m^{\prime}\!>$, 
respectively $\zz\!<\!(E-CA^{-1}B^T)m+(F-CA^{-1}C^T)m^{\prime}\!> \:\leq \zz\!<\!m\!>$, i.e.

\vspace{-5mm}
\begin{eqnarray}\label{lagrange2}
D & = & BA^{-1}B^T\nonumber\\
F & = & CA^{-1}C^T
\end{eqnarray}

\noindent
(for $\qq$-cobordisms this in particular means that the entries on the left-hand side, a priori in
$\zz\left[\frac{1}{\det A}\right]$, must be in $\zz$), and for $\qq$-cobordisms additionally:

\begin{equation}\label{lagrange3}
BA^{-1}C^T \in  {\cal M}_{g_1\times g_2}(\zz)
\end{equation}

\medskip
We will need the following elementary

\enunt{Lemma 12}{ The signature of a symmetric $2g\times 2g$-matrix 
$\left(\begin{array}{cc}A & -I \\ -I & {\bf 0}\end{array}\right)$ with
integer, respectively real entries is $(g,g)$. 
The determinant of such a matrix is $(-1)^g$.}

\medskip
With these notations, Proposition 3 and (\ref{lagrange2}) imply that
the linking matrix $lk(L_1\cup L_0\cup L_2)$ is:

\begin{equation}\label{lkmatrix}
lk(L_1\cup L_0\cup L_2)=\left(
\begin {array}{cccc}
A & B^T & {\bf 0} & {\bf 0} \\
B & BA^{-1}B^T & -I & {\bf 0} \\
{\bf 0} & -I & DC^{-1}D^T & D \\
{\bf 0} & {\bf 0} & D^T & C
\end{array}
\right)
\end{equation}

\noindent 
where $A=lk(L_1)\in{\cal M}_{|L_1|\times|L_1|}(\zz)$, 
$C=lk(L_2)\in{\cal M}_{|L_2|\times|L_2|}(\zz)$, 
$B=lk(G_1^{\prime},L_1)\in{\cal M}_{g\times|L_1|}(\zz)$, 
$D=lk(G_2,L_2)\in{\cal M}_{g\times|L_2|}(\zz)$, 
$BA^{-1}B^T, DC^{-1}D^T\in{\cal M}_{g\times g}(\zz)$. 
With the same notations:

\medskip
\te{Proposition 13}{ The signature of the matrix (\ref{lkmatrix}) is 
$(\sigma_+^1+\sigma_+^2+g,\sigma_-^1+\sigma_-^2+g)$, where $(\sigma_+^1,\sigma_-^1)$,
respectively $(\sigma_+^2,\sigma_-^2)$ is the signature of $lk(L_1)$, respectively
$lk(L_2)$. Also the following holds:
\begin{equation}\label{detgluing}
\det(lk(L_1\cup L_0\cup L_2)) = (-1)^g\cdot\det(lk(L_1))\cdot\det(lk(L_2))
\end{equation}
}

\vspace{-5mm}
\dem{By the classification of quadratic forms with real coefficients there exist matrices
$X\in SL(|L_1|,\rr)$, $Y\in SL(|L_2|,\rr)$ such that $XAX^T$ and $YCY^T$ are diagonal. 
Since $A$ and $C$ are symmetric, so are $BA^{-1}B^T$ and $DC^{-1}D^T$, 
hence there exist $P,Q\in SL(g,\rr)$ such that $PBA^{-1}B^TP^T$ and $QDC^{-1}D^TQ^T$
are diagonal. Then:

$$\left(
\begin {array}{cccc}
X & {\bf 0} & {\bf 0} & {\bf 0} \\
{\bf 0} & P & {\bf 0} & {\bf 0} \\
{\bf 0} & {\bf 0} & Q & {\bf 0} \\
{\bf 0} & {\bf 0} & {\bf 0} & Y
\end{array}
\right)
\left(
\begin {array}{cccc}
A & B^T & {\bf 0} & {\bf 0} \\
B & BA^{-1}B^T & -I & {\bf 0} \\
{\bf 0} & -I & DC^{-1}D^T & D \\
{\bf 0} & {\bf 0} & D^T & C
\end{array}
\right)
\left(
\begin {array}{cccc}
X^T & {\bf 0} & {\bf 0} & {\bf 0} \\
{\bf 0} & P^T & {\bf 0} & {\bf 0} \\
{\bf 0} & {\bf 0} & Q^T & {\bf 0} \\
{\bf 0} & {\bf 0} & {\bf 0} & Y^T
\end{array}
\right)
=$$

$$
=\left(
\begin {array}{cccc}
XAX^T & XB^TP^T & {\bf 0} & {\bf 0} \\
PBX^T & PBA^{-1}B^TP^T & -PQ^T & {\bf 0} \\
{\bf 0} & -QP^T & QDC^{-1}D^TQ^T & QDY^T \\
{\bf 0} & {\bf 0} & YD^TQ^T & YCY^T
\end{array}
\right)
=
\left(
\begin {array}{cccc}
{\cal D}_1 & F_1^T & {\bf 0} & {\bf 0} \\
F_1 & {\cal D}_2 & {\cal E}^T & {\bf 0} \\
{\bf 0} & {\cal E} & {\cal D}_3 & F_2 \\
{\bf 0} & {\bf 0} & F_2^T & {\cal D}_4
\end{array}
\right)
$$

\noindent
where ${\cal D}_i$ are diagonal matrices, ${\cal D}_1$ and ${\cal D}_4$ necessarily invertible, 
and ${\cal E}\in SL(g,\rr)$. Therefore:

$$\left(
\begin {array}{cccc}
I & {\bf 0} & {\bf 0} & {\bf 0} \\
-F_1{\cal D}_1^{-1} & I & {\bf 0} & {\bf 0} \\
{\bf 0} & {\bf 0} & I & -F_2{\cal D}_4^{-1} \\
{\bf 0} & {\bf 0} & {\bf 0} & I
\end{array}
\right)
\left(
\begin {array}{cccc}
{\cal D}_1 & F_1^T & {\bf 0} & {\bf 0} \\
F_1 & {\cal D}_2 & {\cal E}^T & {\bf 0} \\
{\bf 0} & {\cal E} & {\cal D}_3 & F_2 \\
{\bf 0} & {\bf 0} & F_2^T & {\cal D}_4
\end{array}
\right)
\left(
\begin {array}{cccc}
I & -({\cal D}_1^{-1})^TF_1^T & {\bf 0} & {\bf 0} \\
{\bf 0} & I & {\bf 0} & {\bf 0} \\
{\bf 0} & {\bf 0} & I &  {\bf 0}\\
{\bf 0} & {\bf 0} & -({\cal D}_4^{-1})^TF_2^T & I
\end{array}
\right)
=$$

$$
=\left(
\begin {array}{cccc}
{\cal D}_1 & {\bf 0} & {\bf 0} & {\bf 0} \\
{\bf 0} & {\cal D}_2-F_1{\cal D}_1^{-1}F_1^T & {\cal E}^T & {\bf 0} \\
{\bf 0} & {\cal E} & {\cal D}_3-F_2{\cal D}_4^{-1}F_2^T & {\bf 0} \\
{\bf 0} & {\bf 0} & {\bf 0} & {\cal D}_4
\end{array}
\right)
$$

But ${\cal D}_2-F_1{\cal D}_1^{-1}F_1^T
=PBA^{-1}B^TP^T-PBX^T\cdot(X^T)^{-1}A^{-1}X^{-1}\cdot XB^TP^T
={\bf 0}$,
and similarly ${\cal D}_3-F_2{\cal D}_4^{-1}F_2^T
={\bf 0}$. Also observe that

$$
\left(
\begin {array}{cc}
I&{\bf 0}\\
{\bf 0}&-{\cal E}^{-1}
\end{array}
\right)
\left(
\begin {array}{cc}
{\bf 0}&{\cal E}^T\\
{\cal E}&{\bf 0}
\end{array}
\right)
\left(
\begin {array}{cc}
I&{\bf 0}\\
{\bf 0}&-({\cal E}^{-1})^T
\end{array}
\right)
=
\left(
\begin {array}{cc}
{\bf 0}&-I\\
-I&{\bf 0}
\end{array}
\right)
$$

The later matrix, by Lemma 12, has signature $(g,g)$ and determinant $(-1)^g$.
Therefore the signature of the matrix (\ref{lkmatrix}) is 
$(\sigma_+^1+\sigma_+^2+g,\sigma_-^1+\sigma_-^2+g)$.
Since all conjugations above where by matrices of determinant $\pm 1$, the determinant
of the original matrix $(\ref{lkmatrix})$ remained unchanged, 
i.e. we also have the relation $(\ref{detgluing})$.
}

\medskip
As immediate consequences we obtain the two central results of this paper:

\smallskip
\te{Corollary 14}{For semi-Lagrangian $\qq$-cobordisms the integer $s(M,M_1,M_2)$ is always equal to $0$.}

\smallskip
\te{Corollary 15}{For semi-Lagrangian $\qq$-cobordisms $(M,f_1,f_2)$ the cardinality of 
$H_1(\widehat{M},\zz)$ is multiplicative with respect to the composition of cobordisms.}

\medskip
\noindent
{\it Remark.}
For the category $({\mathfrak C}_0,\bullet)$ the same properties, with the integer 
$s:=\sigma^1_++\sigma^2_+-\sigma_+$, are obvious.

\subsection{Consequences}

Set $\T(f|_{\Sigma})=id_{\T(g)}$ for any homeomorphism $f$ of the 
parametrized surfaces. Then $\T$ is a covariant functor, and the naturality axiom (A1) is obvious.

\medskip
\noindent
{\it Basic example.}
Let $N$ be an arbitrary integer, and $\tau(L)\in K^{\ast}$ an invariant of framed links, also invariant under
the second Kirby move, and under changing orientation of link components. It is obvious that, if the linking 
matrix of $L$ is non-singular and has signature $(\sigma_+,\sigma_-)$, and $M\equiv S^3_L$,
then $\tau(M)=|H_1(M,\zz)|\cdot\frac{\tau(L)}{\tau(\oo^{+1})^{\sigma_+}\tau(\oo^{-1})^{\sigma_-}}$
is an invariant of rational homology 3-spheres. 

Assume the link invariant has two additional properties:
\begin{eqnarray}
\tau(\overline{L}) & = & \overline{\tau(L)}\\
\tau(L_1\cup L_0\cup L_2) & = & \tau(L_1)\tau(L_2)\tau(\oo^{+1})^g\tau(\oo^{-1})^g\label{knotrel}
\end{eqnarray}
for any link $L$, $g\in\nn$, and link $L_1\cup L_0\cup L_2$ obtained from $L_1$ and $L_2$ as follows:
add arbitrary $g$ components to $L_i$, and join the two links along the $g$ additional
componenets, preliminary inserting the tangle $T_g$ (figure 4b) in between.
The second property generalizes multiplicativity under connected sum.
The firs property implies that $\tau(-M)=\overline{\tau(M)}$.

Let $\T(g)$, $g\geq 1$ be the $K$-vector space freely generated by homeomorphism classes of
semi-Lagrangian $\qq$-cobordims of the form $(P,s,h)$, 
$s\cup h:\Sigma_0\cup \Sigma_g\stackrel{\approx}{\rightarrow}\partial M$.
The relation $\overline{[P,s,h]}=[P,s,h]$ defines a natural antimorphism on $\T(g)$.
Defines $\tau(M,f,f^{\prime})$ to be the $K$-homomorphism sending
$[P,s,h]\in\T(g)$ to
$\tau(\widehat{M})\cdot[M\cup_{f\circ h^{-1}}P,s,f^{\prime}]\in\T(g^{\prime})$.

(A2) Observe that $\tau(\sigma_g\times[0,1],p_1,p_2)$ sends
$[P,s,h]$ to $\tau(S^3)\cdot[P,s,h]=[P,s,h]$;

(A3) If $M_1$ and $M_2$ are two cobordisms from $\mathfrak Q$, then $\tau(M_2\cup M_1)$ sends $[P,s,h]$ to
$|H_1(M_2\cup M_1)|^N\cdot\frac{\tau(L_1\cup L_0\cup L_2)}{\tau(\oo^{+1})^{\sigma_+}\tau(\oo^{-1})^{\sigma_-}}\cdot
[M_2\cup M_1\cup P,s,f_2^{\prime}]$, while
$\tau(M_2)\circ\tau(M_1)$ sends it to
$|H_1(M_2)|^N\cdot |H_1(M_1)|^N\cdot
\frac{\tau(L_1)\tau(L_2)}{\tau(\oo^{+1})^{\sigma_+^1+\sigma_+^2}\tau(\oo^{-1})^{\sigma_-^1+\sigma_-^2}}\cdot
[M_2\cup M_1\cup P,s,f_2^{\prime}]$. 
These are equal by (\ref{knotrel}) and Corollaries 14 and 15;

(A4) For $[Q,r,j]^{\ast}\in\T(g)$ define the dual map $[Q,r,j]^{\ast}$ to send $[P,s,h]$
to $\tau(\widehat{Q})\cdot\tau(\widehat{-Q\cup_{j\circ h^{-1}}P,s,-r})$.
Note that if one represents $(M,f,f^{\prime})$ by a triplet (Proposition 1), takes a generic projection on $\rr^2$,
change all crossing to their opposites, and apply $\kappa$, one obtains the cobordism $(-M,-f^{\prime},-f)$. 
Therefore defining $\overline{\tau(M,f,f^{\prime})}=\tau(-M,-f^{\prime},-f)$ we obtain for the case
when the domain of $f$ is $\Sigma_0$ the composition $proj(\ast(\overline{\cdot}))$, 
where $proj:\T(0)\rightarrow K$, $proj[P,s,r]=\tau(\widehat{P})$.

If we replace $\T(0)$ by its image through $proj$,
the above data defines a non-degenerate anomaly-free TQFT on $\mathfrak Q$. 

\medskip
Of cause, given a specific $\tau$, one should define $\T(g)$, $\tau(M,f,f^{\prime})$, the composition and conjugation
in a way directly related to $\tau$ of links. Non-degeneracy should also be addressed in $\tau$-spceific terms.
In [\ref{CL}] we use the results of this paper to construct a TQFT for the Le-Murakami-Ohtsuki invariant of $\qq$HS,
as well as for its degree truncations. 

\medskip\noindent
{\it Remark.}
Throughout this paper we have taken {\it oriented} chain graphs. This is motivated by the fact that for our main example
(LMO) absence of orientation would induce unjustified complications. However for the correspondence in Proposition 1
here it sufices to consider banded (rather than framed) graphs with a specification of $g$ meridional disks for
its circle components (which would replace our notion of "horizontal line").

\medskip
The TQFTs based on $\mathfrak Q$, when restricted to $\mathfrak L$ produce linear representations
${\cal L}_g\rightarrow GL_{K}(\T(g))$.
The group ${\cal L}_g$ has not been studied before, no set of relations\footnote{A set of generators can
be easely obtained by taking products of generators of $\T_g$ and of $GL(2,\zz)$} is known. 
Even, if we restrict to Torelli group, although a finite set of generators for 
${\cal T}_g$ is well-known [\ref{johnson}], existance of a finite presentation is an open problem.


%
%
%
%
%
%
%
%
%
%
%
%
%
%
%
%
%

\section*{References}
\sloppy

\begin{bib} 
\item\label{at88} 
\rm M.Atiyah, \it Topological Quantum Field Theories, \rm
Publications Math\'{e}matiques IHES {\bf 68} \rm , 175-186 (1988)

\item\label{BHMV}
\rm C.Blanchet, H.Habegger, G.Masbaum, P.Vogel, \it Topological
quantum field theories derived from the Kauffman bracket, \rm
Topology {\bf 34}, no.4, 883-927 (1995)

\item\label{CL}
\rm D.Cheptea, T.Le, \it A TQFT associated to LMO invariant of three-dimensional manifolds,
\rm math.GT.0508220 v2

\item\label{deloup}
\rm F.Deloup, \it An explicit construction of an abelian topological 
quantum field theory in dimension 3, \rm Topology and Its Applications 
{\bf 127}, no 1-2, 199-211 (2003) 

\item\label{FM}
\rm A.Fomenko, S.Matveev, \it Algorithmic and computer methods for 
three-manifolds, \rm Kluwer Academic Publishers (1997) 

\item\label{gervais}
\rm S.Gervais, \it A finite presentation of the mapping class
group of a punctured surface, \rm Topology {\bf 40}, no.4, 703-725 
(2001)

\item\label{GS}
\rm R.E.Gompf, A.I.Stipsicz, \it 4-manifolds and Kirby calculus, 
\rm Graduate Studies in Mathematics {\bf 20}, AMS (1999)

\item\label{habiro}
\rm K.Habiro, \it Claspers and finite-type invariants of links, \rm  Geometry
and Topology {\bf 4}, 1-83 (2000) 

\item\label{johnson}
\rm D.Johnson, \it The structure of the Torelli group I: A finite set of
generators for ${\cal T}$, \rm Annals of Mathematics {\bf 118}, 423-442 (1983)

\item\label{LMO}
\rm T.T.Q.Le, J.Murakami, T.Ohtsuki, \it On a universal perturbative
invariant of 3-manifolds, \rm Topology {\bf 37}, no.3, 539-574 (1998)

\item\label{matveev}
\rm S.V.Matveev, \it Generalized surgery of three-dimensional 
manifolds and representations of homology spheres, \rm Matematicheskie
Zametki {\bf 42}, no.2, 268-278 (1986)

\item\label{milnor}
\rm J.Milnor, \it Lectures on the h-cobordism theorem, 
\rm Princeton University Press (1965)

\item\label{morita}
\rm S.Morita, \it Casson's invariant for homology 3-spheres and 
characteristic classes of surface bundles I, \rm Topology {\bf 28}, no.3, 
305-323 (1989)

\item\label{MO}
\rm J.Murakami, T.Ohtsuki, \it Topological Quantum Field Theory
for the Universal Quantum Invariant, \rm Commun. Math. Phys. {\bf
188}, 501-520 (1997)

\item\label{turaev}
\rm V.Turaev, \it Quantum Invariants of Knots and 3-Manifolds, \rm
Walter de Gruyter (1994)

\item\label{vogel}
\rm P.Vogel, \it Invariants de type fini, \rm en ``Nouveaux
Invariants en G{\' e}om{\' e}trie et en Topologie'', publi{\' e}
par D. Bennequin, M. Audin, J. Morgan, P. Vogel, Panoramas et
Synth{\` e}ses {\bf 11}, Soci{\' e}t{\' e} Math{\' e}matique de
France , 99-128 (2001)

\end{bib}

\begin{flushleft}
Dorin Cheptea\newline
\sc UFR de Math\'ematique et d'Informatique, Universit\'e Louis Pasteur,\newline
7, rue Ren\'e Descartes, 67084, Strasbourg, France\newline
and\newline
\sc Institute of Mathematics, P.O.Box 1-764, Bucharest, 70700, Romania \newline
e-mail: \tt cheptea@math.u-strasbg.fr

\medskip\rm
Thang T Q Le\newline
\sc School of Mathematics, Georgia Institute of Technology,\newline
Atlanta, GA 30332-0160, USA \newline
e-mail: \tt letu@math.gatech.edu
\end{flushleft}

\end{document}